\documentclass[EJP]{ejpecp} % replace ECP by EJP if needed.
\usepackage{mathdots,relsize}
\SHORTTITLE{Random Power Towers}

\TITLE{Infinite Random Power Towers} % \thanks is optional. Insert line breaks with \\

\AUTHORS{%
  Mark~Dalthorp\footnote{Cornell University, United States of America.
    \EMAIL{md944@cornell.edu}}}%AUTHORS
 % Separate items with ;

\AMSSUBJ{60J05} % Edit. Separate items with ;
\SUBMITTED{October 24, 2021} % Edit.
\ACCEPTED{January 3, 2024} % Edit.

\VOLUME{0}
\YEAR{2023}
\PAPERNUM{0}
\DOI{10.1214/YY-TN}

\ABSTRACT{We prove a probabilistic generalization of the classic result that infinite power towers, $c^{c^{\iddots}}$, converge if and only if $c\in[e^{-e},e^{1/e}]$. Given an i.i.d. sequence $\{A_i\}_{i\in\nn}$, we find that convergence of the power tower $A_1^{A_2^{\iddots}}$ is determined by the bounds of $A_1$'s support, $a=\inf(\supp(A_1))$ and $b=\sup(\supp(A_1))$. When $b\in[e^{-e},e^{1/e}]$, $a<1<b$, or $a=0$, the power tower converges almost surely. When $b<e^{-e}$, we define a special function $B$ such that almost sure convergence is equivalent to $a<B(b)$. Only in the case when $a=1$ and $b>e^{1/e}$ are the values of $a$ and $b$ insufficient to determine convergence. We show a rather complicated necessary and sufficient condition for convergence when $a=1$ and $b$ is finite.
	
We also briefly discuss the relationship between the distribution of $A_1$ and the corresponding power tower $T=A_1^{A_2^{\iddots}}$. For example, when $T\sim\mathrm{Unif}[0,1]$, then the corresponding distribution of $A_1$ is given by $UV$ where $U,V\sim\mathrm{Unif}[0,1]$ are independent. We generalize this example by showing that for $U\sim\mathrm{Unif}[\alpha,\beta]$ and $r\in\rr$, there exists an i.i.d. sequence $\{A_i\}_{i\in\nn}$ such that $U^r \eqindist A_1^{A_2^{\iddots}}$ if and only if $r\in[0, \frac1{1+\log \beta}]$.}

\newcommand{\nn}{\mathbb{N}}
\newcommand{\rr}{\mathbb{R}}
\newcommand{\zz}{\mathbb{Z}}
\newcommand{\cc}{\mathbb{C}}
\newcommand{\ilim}[1]{\lim\limits_{#1\rightarrow\infty}}
\newcommand{\dd}[1]{\frac{d}{d #1}}

\newcommand{\ee}{\mathbb{E}}
\newcommand{\supp}{\mathrm{supp}}
\newcommand{\floor}[1]{\lfloor #1 \rfloor}
\newcommand{\abs}[1]{\left| #1 \right|}
\newcommand{\toindist}{\stackrel{d}{\rightarrow}}
\newcommand{\toinp}{\stackrel{P}{\rightarrow}}
\newcommand{\eqindist}{\stackrel{d}{=}}
\newcommand{\at}[1]{\textrm{\footnotesize \textnormal {AT}}_{#1}}
\newcommand{\ate}{\at{\mathrm{even}}}
\newcommand{\ato}{\at{\mathrm{odd}}}
\newcommand{\et}{\mathop{\mathlarger{\mathrm E}}}

\newtheorem{theoremcase}{Theorem \ref{th:main} case}

\begin{document}

%%%%%%%%%%%%%%%%%%%%%%%%%%%%%%%%%%%%%%%%%%%%%%%%%%%%%%%%%%%%%%%%%%%
%%                                                               %%
%% No need for \maketitle.                                       %%
%%                                                               %%
%%%%%%%%%%%%%%%%%%%%%%%%%%%%%%%%%%%%%%%%%%%%%%%%%%%%%%%%%%%%%%%%%%%

%%%%%%%%%%%%%%%%%%%%%%%%%%%%%%%%%%%%%%%%%%%%%%%%%%%%%%%%%%%%%%%%%%%
%%                                                               %%
%% Please replace what follows by the body of your article       %%
%% (up to the bibliography):                                     %%
%%                                                               %%
%%%%%%%%%%%%%%%%%%%%%%%%%%%%%%%%%%%%%%%%%%%%%%%%%%%%%%%%%%%%%%%%%%%

\section{Introduction and Main Result}

\subsection{Background}
In this paper we investigate infinite random power towers, defined as the limit of the sequence $$
A_1 ,\, (A_1)^{A_2},\,(A_1)^{(A_2)^{A_3}},\dots
$$
where $(A_i)_{i\in\nn}$ is an i.i.d. sequence of positive real random variables. We determine necessary and sufficient conditions for convergence of the power tower which apply unless the $A_1$ has support bounded below by $1$ and unbounded above. We also look at some specific cases where the distribution of the limit may be written in closed form. 

This question comes as a natural probabilistic extension of the long history of study of power towers, which we outline briefly below. Power towers have usually been studied as a topic of inherent interest, but they do occasionally find applications in other topics (see for example \cite{article:9},\cite{Dzurina}, or \cite{knoebel}). They also provide a simple example of a non-trivial discrete dynamical system. Similarly, the random power tower provides a simple example of a random dynamical system. In this paper we will see that the fact that exponential functions are monotone allows for some shortcuts not possible in full generality. 

The simplest infinite power towers, taking the form \[
c^{c^{c^{\iddots}}}
\]
were first studied by Euler in 1778 \cite{ARTICLE:4} and independently a few decades later by Eisenstein in 1844 \cite{ARTICLE:2}. The infinite expression is not \emph{a priori} well-defined, but must be defined as the limit of the sequence \[
c, \hspace{4pt} c^c, \hspace{4pt} c^{c^{c}},\hspace{4pt}\dots
\]
On first exposure to power towers, it may be surprising that it is possible for them to converge at all when $c>1$. After all, exponentiation famously diverges to infinity rapidly, one would expect the iterated sequence to diverge even more rapidly. However, for nonnegative $c$, the infinite power tower converges if and only if $c\in[e^{-e},e^{1/e}]$ \cite{ARTICLE:4}. Even if the existence of this interval of convergence isn't a surprise, the endpoints seem far nicer than they have a right to be.

\begin{figure}[t]
	\centering
	\includegraphics[scale=0.54]{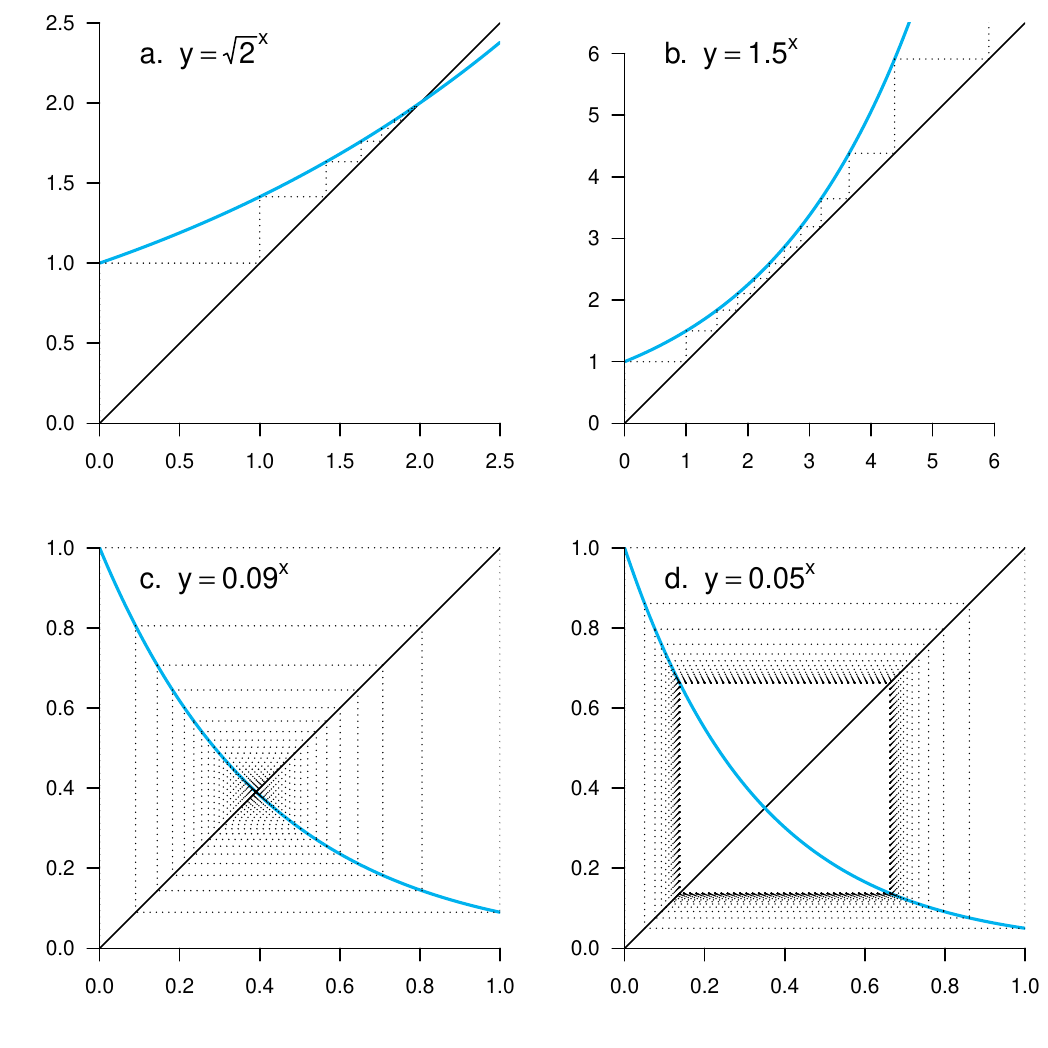}\caption{Iterating exponential functions showing the four possible behaviors: a. Increasing convergence to a fixed point when $c\in[1,e^{\frac1e}]$; b. Increasing to infinity when $c>e^{\frac1e}$; c. Alternating convergence to a fixed point when $c\in[e^{-e},1]$; d. Oscillation in the limit when $c<e^{-e}$.}\label{fi:sqrtower}
\end{figure}

The sequence $c,c^c,c^{c^c},\dots$ is naturally viewed as repeated iterates of the function $x\to c^x$ starting at $x=1$, which turns the problem into a discrete dynamical system. The theory of dynamical systems is complex and deep; in general their behavior is difficult or impossible to predict. However, this particular system is easy to get a handle on. For example, the iteration of $x\to(\sqrt 2)^x$ forms a sequence that increases to the first fixed point of this function, i.e. $2$ (Figure \ref{fi:sqrtower}).

\begin{figure}[t]
	\centering
	\includegraphics[scale=0.5]{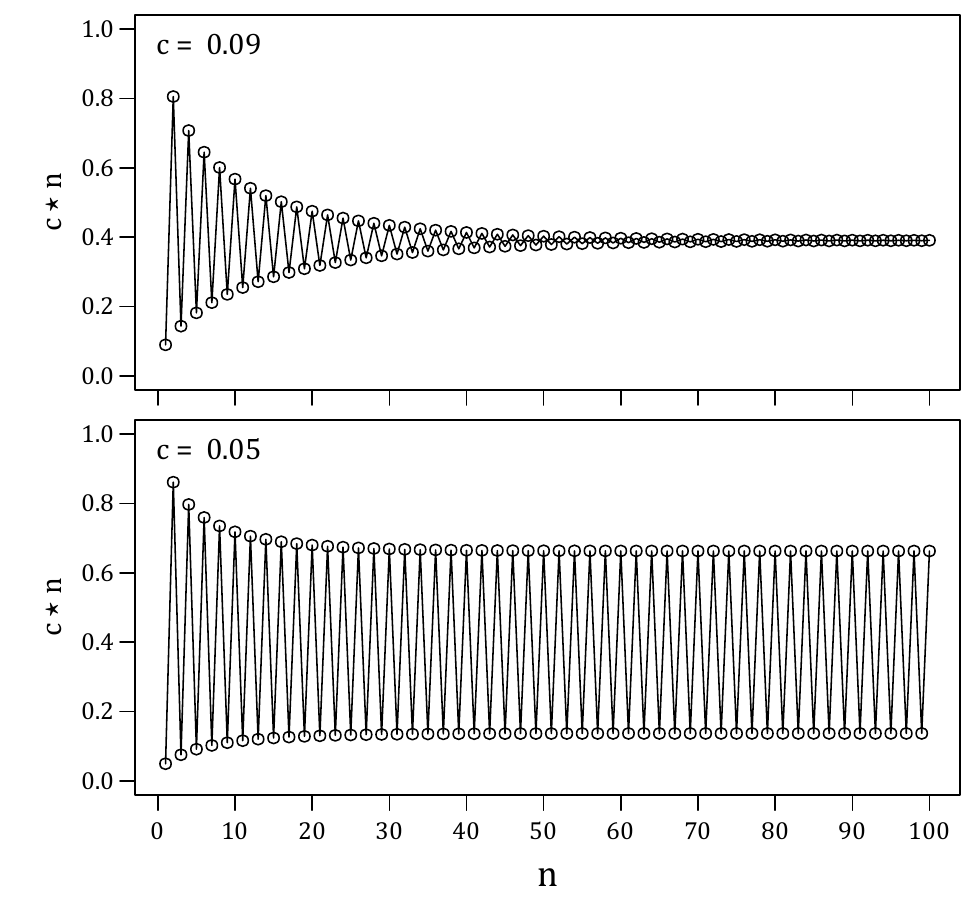}\caption{Limiting behavior of $c,c^c,c^{c^c},\dots$ for $c$ on either side of $e^{-e}\approx 0.06$}\label{fi:oscillation}
\end{figure}

In general, iterating an increasing function creates sequence that either converges to a fixed point or diverges to infinity. Sequences arising from iterated decreasing functions will either converge to a fixed point, or oscillate in the limit. The convergence results about $c^{c^{\dots}}$ become perfectly natural in light of these facts: The sequence $c, c^c, c^{c^c}, \dots,$ converges in the limit when $c\in[e^{-e}, e^{\frac1e}]$. For larger $c$, the sequence diverges to infinity. For $c\in (0,e^{-e})$, it oscillates, with the even- and odd-indexed subsequences converging to distinct limits (Figure \ref{fi:oscillation}). We will see similar behavior in the random case as discussed below.

This is a special case of the fully general question of convergence of $$c_1,\, (c_1)^{c_2},\,(c_1)^{(c_2)^{c_3}},\dots$$ for an arbitrary sequence $c_i$. In this paper, we will restrict our attention to the case of positive real sequences (the complex case was investigated by Thron \cite{Thron}). To our knowledge, this general question was first posed in the literature by D.F. Barrow in 1936 \cite{Barrow}, who proved that this sequence converges if $c_i\in[e^{-e},e^{1/e}]$, among other interesting results. Barrow also introduced the notation $$
\et_{k=1}^n c_k = c_1^{c_2^{{\cdots}^{c_n}}} 
$$
and $$
\et_{k=1}^\infty c_k = \ilim n \et_{k=1}^n c_k
$$
by analogy to summation and product notation, which we will follow here. Baker and Rippon \cite{ARTICLE:7} gave necessary and sufficient conditions for convergence of an alternating power tower with positive real coefficients, i.e. $c_k$ alternates between $a$ and $b$, with $a,b\in(0,\infty)$. Later \cite{ARTICLE:8}, they found conditions for convergence and divergence when $c_k$ is periodic with other periods. Bachman (1995) \cite{Bachman} proves an essentially tight bound on convergence in the event that $c_k$ converges to $e^{1/e}$ from above, based on a note without proof found in Ramanujan's notebook.

In this paper, we investigate the probabilistic question: If $\{A_i\}_{i\in \nn}$ is an i.i.d. sequence of random variables, when does $\et\limits_{i=1}^\infty A_i$ converge? In our main result, Theorem \ref{th:main}, we prove necessary and sufficient conditions for almost sure convergence under the assumption that the support of $A_1$ is a bounded subset of $[0,\infty)$. Surprisingly, almost sure convergence depends only on the upper and lower bounds of the support of $A_1$, except when the lower bound is $1$ and the upper bound is larger than $e^{1/e}$. Along the way in our proof, we will derive a closed form for a function whose existence was proven by Baker and Rippon \cite{ARTICLE:7}, which defines the boundary between convergence and divergence of alternating power towers $a^{b^{a^{b^{\cdots}}}}$ for $a,b\in(0,1)$. 

We end by asking which distributions arise as a limit of a sequence of random power towers: Given a random variable $X$, can we find an i.i.d. sequence $A_i$ such that $X = \et_{i=1}^\infty A_i$? If so, we will say that $X$ has the tower property, and the distribution of $A_i$ is its inverse tower distribution. In Theorem \ref{th:inverse_problem} we prove that if $U\sim\mathrm{Unif}[\alpha,\beta]$ and $r\in\rr$, then $U^r$ has the tower property if and only if $1\in[\alpha,\beta]$ and $r\in[0,\frac1{1+\log \beta}]$.

As in the deterministic case, it is natural to think of random power towers as composing random maps $x\to (A_n)^x$. Composition sequences of random maps has also been studied extensively, for example, see the review by Diaconis and Freedman \cite{ARTICLE:1}, who showed a probabilistic analogue of the Banach fixed point theorem, which may be applied in many diverse instances of random iterated functions.  Stated informally, if $f_n$ is a random i.i.d. sequence of functions from some metric space $S$ to itself, and $f_n$ is ``on average" a contraction, then \[
f_1,\hspace{4pt}  f_1\circ f_2,\hspace{4pt}  f_1\circ f_2\circ f_3,\hspace{4pt}  \dots
\] almost surely converges to a constant function in the limit. This is a powerful result, but it is fairly limited for our question. The best that this can give us is
\begin{corollary}[Of Proposition 5.1 in \cite{ARTICLE:1}]\label{cor:intro}
	If $\{A_n\}_{n\in\infty}$ is an i.i.d. sequence of random variables on $[0,e^{\frac1e}]$ such that \[
	\ee\left[\log\left(\max(-\log A_1, (A_1)^e\log A_1)\right)\right] < 0
	\]
	then the infinite random power tower $\et\limits_{i=1}^\infty A_i$ converges almost surely.
\end{corollary}
Note that the expression in the expectation is simply the log of Lipschitz constant of the map $x\rightarrow (A_1)^x$ on the interval $[0,e]$. The condition of the log Lipschitz constant having negative expectation was called ``super-contracting" by Steinsaltz \cite{ARTICLE:5}, who also studied infinite iterated function systems, and gave convergence results with weaker assumptions than Diaconis and Freedman's theorem. If we allow the base to exceed $e^{\frac1e}$, then their result is not much help, as $x\to a^x$ is not a Lipschitz map of any subinterval of $[0,\infty)$ to itself. For bases close to $0$, we have the difficulty that $\log|\log a|$ becomes very large, so one might expect the random power tower to have an oscillating limit if most of its weight is near $0$, similar to the behavior of $a^{a^{\iddots}}$ when $a<e^{-e}$.

There are in fact many distributions for $A_1$ that will give rise to an almost-surely-convergent power tower that do not satisfy the ``contracting on average" condition, and indeed that do not satisfy any condition based on Lipschitz constants.

\subsection{Notation}
As mentioned above, we will be following Barrow's power tower notation with $\et\limits_{i=1}^n c_i = (c_1)^{(c_2)^{{\dots}^{c_n}}}$ and $\et\limits_{i=1}^\infty c_i = (c_1)^{(c_2)^{{\dots}}}$. We will also find it convenient to add to these \[
\et\limits_{i=1}^n c_i ; x = (c_1)^{(c_2)^{{\dots}^{(c_n)^x}}}.
\]
Similar to exponentiation, we define this notation to associate to the right, i.e. $$\et\limits_{i\in A} a_i;\et\limits_{j\in B}b_j;\et\limits_{k\in C} c_k = \et\limits_{i\in A} a_i;\left(\et\limits_{j\in B}b_j;\left(\et\limits_{k\in C} c_k\right)\right)
$$

Also, by analogy to defining the empty product as $1$, we define the empty power tower $\et_{i=n+1}^n c_i = 1$.

For the sake of concision, we will want to define
\[c\star n=\et_{i=1}^n c = \underbrace{c^{c^{\iddots^c}}}_{n \mbox{ times}}\] This operation is commonly called ``tetration" and is often denoted ${^n}c$ or $c\uparrow \uparrow n$, but we find this notation easier to read and pronounce (``$c$ star $n$"). We also define the iterated logarithm function \[
\log^\star (x) = \inf\{n\in\nn : e\star n \ge x \}
\]
which is a right-inverse of $n \rightarrow e\star n$, i.e. $\log^\star (e\star n) = n$ for $n\in \nn$. Additionally, $e\star \log^\star (x) = \floor{x}$. Just as $e\star n$ diverges to infinity extremely rapidly, so $\log^\star (x)$ diverges to infinity extremely slowly. The iterated logarithm function is used in study of algorithms, and this notation is fairly standard there \cite{BOOK:1}.

Also, let $W(z)$ be the principal branch of the Lambert-$W$ function, defined to be the inverse function of $z e^z$ on $\cc$. As a real-valued function, $W(z)$ has domain $[-\frac1e,\infty]$ and range $[-1,\infty]$. See \cite{ARTICLE:6} for an in-depth introduction to this function.

Throughout, we will always have $\{A_n\}_{n\in \nn}$ representing an i.i.d. sequence on $(0,\infty)$, and $T_n = \et\limits_{i=1}^n A_i$, and $a=\inf(\supp(A_1))$ and $b=\sup(\supp(A_1))$. 

In this paper, any unqualified equations or inequalities involving random variables may be taken either absolutely or almost surely, depending on whether the involved random variables are surely between the endpoints of their support, or only almost surely in that range.
\subsection{The main theorem}
The following theorem gives necessary and sufficient conditions for almost-sure convergence of the random power tower $T_n$ provided that the support of $A_1$ is not an unbounded subset of $[1,\infty]$:
\begin{theorem}\label{th:main} We have four different cases, depending on the values of $a$ and $b$:\begin{enumerate}
		\item If $[a,b]\subseteq[1,e^{1/e}]$, then $T_n$ converges almost surely.
		\item If $1\le a$ and $ b \in (e^{1/e},\infty)$:
		Then $T_n$ converges a.s. if \[\ee\left[\inf\{n: A_n \le e^{\frac 1 { e\star n}}\} \right]< \infty\]
		and $T_n$ diverges to infinity a.s. if \[\ee\left[\inf\{n: A_n \le e^{\frac 1 { e\star n}}\} \right]= \infty\]
		\item Let $B:[0,1]\rightarrow[0,1]$ be defined by \[
		B(x) =\begin{cases} x& x\ge e^{-e}\\\exp\left(W\left(\frac1{\ln x}\right)\exp\left(-\frac1{W\left(\frac1{\ln x}\right)}\right)\right) & x<e^{-e}\end{cases}
		\]
		If $b\le1$ and $a<B(b)$, then $T_n$ converges almost surely. Conversely, if $b<e^{-e}$ and $a\ge B(b)$, then $T_n$ diverges by oscillation, with the even and odd subsequences converging to distinct limits.
		\item
		If $a<1<b\le\infty$, then $T_n$ converges almost surely.
	\end{enumerate}
\end{theorem}
For the purposes of checking given a distribution for $A_1$ whether $T_n$ converges, the following corollaries are easier to use (especially since we can do without the rather technical condition of case 2):
\begin{corollary}[Convergence conditions]\label{cor:conv} If any of the following holds, then $T_n$ converges almost surely:\begin{enumerate}
		\item $a=0$
		\item $b\in[e^{-e},e^{\frac1e}]$
		\item $a<1\le b$
		\item $b<e^{-e}$ and \[
		a < \exp\left(W\left(\frac1{\ln b}\right)\exp\left(-\frac1{W\left(\frac1{\ln b}\right)}\right)\right)
		\]
		\item $a=1$ and $b$ is finite and  \[
		\liminf\limits_{t\rightarrow 1^+} \left[P(A_1 \le t)\log^\star\left(\frac1{t-1}\right)\right] > 1
		\]
	\end{enumerate}
\end{corollary}
\begin{corollary}[Divergence conditions] If any of the following holds, then $T_n$ diverges almost surely:\begin{enumerate}
		\item $a>1$ and $b>e^{\frac1e}$
		\item $b<e^{-e}$ and \[
		a \ge \exp\left(W\left(\frac1{\ln b}\right)\exp\left(-\frac1{W\left(\frac1{\ln b}\right)}\right)\right)
		\]
		\item $a=1$ and $b>e^{\frac1e}$ and \[
		\ee \left[\log^\star\right(\frac1{A_1 - 1}\left)\right] < \infty
		\]
            \item $a=1$ and $b>e^{\frac1e}$ and  \[
		\limsup\limits_{t\rightarrow 1^+} \left[P(A_1 \le t)\log^\star\left(\frac1{t-1}\right)\right] < 1
		\]
	\end{enumerate}
	
\end{corollary}

Surprisingly, unless $a=1$, the particulars of the distribution of $A_1$ are completely irrelevant except for the bounds of its support. One can also see immediately that this is a much more powerful result than the ``contracting on average" condition from the introduction; indeed $T_n$ converges almost surely when $a=0$ and $b<e^{-e}$, but $x\to (A_n)^x$ is never a contraction.

Even when $a=1$ and $b\in (e^{1/e},\infty)$, the only fact about $A_1$ which matters is how much the distribution of $A_1$ is weighted on neighborhoods around $1$. Interestingly, because $\log^\star\left(\frac1{t-1}\right)$ goes to infinity so slowly as $t$ approaches $1$ from above, for any function $f:[1,\infty)\rightarrow[1,\infty)$ with $f(1)=1$ that is Hölder continuous at $1$ and bounded on $[1,b]$, convergence of $T_n$ is implies convergence of $\et_{k=1}^\infty f(A_k)$. 

There is no reason to expect that any ``nice" distribution for $A_1$ should correspond to a nice distribution of $T$, but we do have the following fortuitous example:

\begin{example}\label{ex1}
	If $U_n$ and $V_n$ are two independent i.i.d. sequences, uniformly distributed on $(0,1)$, then the sequence $\et_{i=1}^n (U_iV_i)$ converges almost surely, and the limit $\et_{i=1}^\infty (U_iV_i)$ is uniform on $(0,1)$. 
\end{example}

Almost sure convergence follows immediately from the first convergence condition of Corollary \ref{cor:conv}, even though this distribution fails to be ``contracting on average".

Proof that its limit is uniform is slightly more involved: If $U,V$ and $W$ are all independent uniformly distributed random variables on $(0,1)$, computation of the elementary triple integral $P((UV)^W \le x)$ shows that $(UV)^W \eqindist W$. As a result of this fact, we have that $\et_{i=1}^{2n}(U_iV_i);U_{2n+1}$ is uniformly distributed on $(0,1)$ for each $n$. We furthermore have
\[\et_{i=1}^{2n+1} (U_iV_i) \le \et_{i=1}^{2n}(U_iV_i);U_{2n+1} \le \et_{i=1}^{2n} (U_iV_i),\]
which follows from the fact that the leftmost quantity equals $\et_{i=1}^{2n}(U_iV_i);(U_{2n+1}V_{2n+1})$, and $\et_{i=1}^{2n}(U_iV_i); x$ is an increasing function of $x$. The leftmost and rightmost terms both converge to the same limit, namely $\et_{i=1}^\infty (U_iV_i)$, hence the middle term also converges the same thing. Since the middle term is uniformly distributed on $(0,1)$ for each $n$, it follows that its limit is as well. 
\section{Proof of Theorem \ref{th:main}}
As we have seen, we cannot use contraction-type of results to prove Theorem \ref{th:main}. However, exponential functions do have a different useful property: They are always monotonic. The proof in all different cases relies heavily on this fact, but it works differently when $a\ge 1$ and when $a<1$, because in the former case, $x\to (A_n)^x$ is always nondecreasing, whereas in the latter case, it can be either nonincreasing or nondecreasing depending on which side of $1$ the base is on.

\subsection{Cases 1 and 2 $(a\ge 1)$}
Importantly, $a\ge 1$ implies that $T_n$ is a non-decreasing sequence because \[
T_{n+1} = A_1^{A_2^{{\iddots}^{A_n^{(A_{n+1})}}}} \ge A_1^{A_2^{{\iddots}^{A_n^{1}}}} = T_n
\]
Our proof strategy for convergence in cases 1 and 2 involves putting an appropriate upper bound on $T_n$ to ensure convergence. For case 1, this is trivial:
\begin{theoremcase} If $[a,b]\subseteq[1,e^{1/e}]$, then $T_n$ converges almost surely.\end{theoremcase}
\begin{proof}
This admits several trivial proofs. It is immediate consequence of Theorem 6 in Barrow \cite{Barrow} or of Corollary \ref{cor:intro} from our introduction. Alternatively, we may prove it directly: By the above, $T_n$ is non-decreasing. It is bounded above by $e$, because for $x\in[1,e^{1/e}]$ and $y\le e$, $x^y \le e$. Therefore $T_n$ converges.\end{proof}

\begin{theoremcase}If $1\le a$ and $ b \in (e^{1/e},\infty)$:
	Then $T_n$ converges a.s. if \[\ee\left[\inf\{n: A_n \le e^{\frac 1 { e\star n}}\} \right]< \infty\]
	and $T_n$ diverges to infinity a.s. if \[\ee\left[\inf\{n: A_n \le e^{\frac 1 { e\star n}}\} \right]= \infty\]
\end{theoremcase}
\emph{Proof.} The proof is in several steps. The basic idea here is that, any time a run $A_i,A_{i+1},\dots,A_{j}$ are all far from $1$, the partial power tower $\et_{k=i}^j A_k$ becomes very large, but if $A_{i-1}$ is extremely close to $1$, it will counteract the run of large terms and $\et_{k=i-1}^j A_k$ may still be close to $1$. Thus if $A_n$ is close enough to $1$ often enough, we may have convergence of $T_n$, but if not, we expect $T_n$ to diverge to infinity.

The quantity $\ee\left[\inf\{n: A_n \le e^{\frac 1 { e\star n}}\} \right]$ actually measures how frequently $A_n$ is sufficiently close to $1$. Because $e^{\frac 1 { e\star n}}$ converges to $1$ rapidly, we have $P(A_n \le e^{\frac 1 { e\star n}})$ goes to $0$ as well (unless $P(A_1=1)>0$), and this expectation being finite means that these probabilities do not go to $0$ too rapidly. For our proof, we substitute this condition for another, slightly easier to compute one, then show equivalence of the two at the end: We start by showing that convergence of $T_n$ is equivalent to convergence of \begin{equation}\label{eq:sum_check}
\sum_{k=1}^\infty \prod_{j=1}^{k-1} P(A_j>b^{\frac1{b\star j}}).
\end{equation}
For convenience, we define $p_j = P(A_n > b^{\frac1{b\star j}})$, which is well-defined because $A_n$ are identically distributed.
	
\subsubsection{Summability of (\ref{eq:sum_check}) implies $T_n$ converges}
Here, we use the convergence of (\ref{eq:sum_check}) to construct a distributional upper bound for $T_n$. For deterministic nondecreasing sequences, proving existence of an upper bound proves convergence. We use the following lemma, which is a probabilistic analogue of this fact, to prove convergence of $T_n$.
\begin{lemma}\label{lem:bounded}
	Let $\{X_n\}_{n\in\nn}$ be an almost surely non-decreasing sequence of random variables. Suppose that there exists an identically distributed (but not necessarily independent) sequence of random variables $\{Y_n\}_{n\in\nn}$ such that $X_n\le Y_n$ almost surely for all $n$. Then $
	\lim\limits_{n\rightarrow\infty}X_n
	$
	converges with a nonzero probability.
\end{lemma}
In the case of $\ee Y_n < \infty$, this result follows from Markov's inequality. The general version is a straightforward modification.

Notice that $T_n$'s convergence is a tail event (Theorem 1 in \cite{Barrow}), and therefore it either has probability $1$ or $0$ by Kolmogorov's 0-1 Law. Hence, it suffices to show that $T_n$ converges with nonzero probability. By Lemma \ref{lem:bounded}, it suffices to construct a sequence of random variables $B_n$, constant in distribution, such that $T_n\le B_n$. This construction will be fairly involved:

Begin by defining a triangular array of random integer-valued random variables $N_{\ell,n}$ for $0\le\ell\le n$. We define $N_{0,n}$ to be i.i.d. and independent from $\{A_i\}_{i\in \zz}$, each having distribution \[
	P(N_{0,n}=m) = \frac{\prod\limits_{j=1}^{m-1} p_j}{\sum\limits_{k=1}^\infty \prod\limits_{j=1}^{k-1} p_j}.
	\]
	This is a well-defined probability distribution by the assumption that $\sum\limits_{k=1}^\infty \prod\limits_{j=1}^{k-1} p_j<\infty$. For $\ell>0$, define $N_{\ell,n}$ recursively by \[
	N_{\ell+1,n} = \begin{cases}1 & \text{if }A_{n-\ell} \le b^{\frac1{b\star N_{\ell,n}}}\\ 1+N_{\ell,n} & \text{otherwise}\end{cases}.
	\]
	Note that $N_{\ell,n}$ becomes $1$ when $A_{n-\ell+1}$ is close to $1$, and otherwise increments by $1$. In this way, $N_{\ell,n}$ is (roughly speaking) a measure of the length of a run of "distant-from-1" values of $A_i$ starting with index $i=n-\ell+1$. The requisite closeness to $1$ is determined by how long the subsequent run of large values is.
	
	We will ultimately use $B_n = b\star N_{n,n}$ as our bound on $T_n$. First, we show that $N_{\ell,n}$ has identical distribution for all $\ell\in \{0,1,\dots,n\}$ by inducting on $\ell$. Suppose \[
	P(N_{\ell,n} = m) = \frac{\prod\limits_{j=1}^{m-1} p_j}{\sum\limits_{k=1}^\infty \prod\limits_{j=1}^{k-1} p_j}
	\]
	for all $m$. This is true by definition for $\ell = 0$. For $\ell >0$ and $m=1$: \begin{eqnarray}
		P(N_{\ell+1,n}=1) &=& P(A_{n-\ell}\le b^{\frac1{b\star N_{\ell,n}}})= \sum_{k=1}^\infty P(N_{\ell,n} = k)P(A_{n-\ell}\le b^{\frac1{b\star N_{\ell,n}}}\mid N_{\ell,n}=k) \nonumber\\
		&=& \sum_{k=1}^\infty P(N_{\ell,n}=k) (1-p_k) =\sum\limits_{k=1}^\infty\frac{ \prod\limits_{j=1}^{k-1}p_j}{\sum\limits_{k=1}^\infty \prod\limits_{j=1}^{k-1} p_j}(1-p_k) = \frac{1}{\sum\limits_{k=1}^\infty \prod\limits_{j=1}^{k-1} p_j} .\nonumber
	\end{eqnarray}
	For $\ell >0$ and $m>1$, we again have $P(N_{\ell+1,n} = m) = \sum_{k=1}^\infty P(N_{\ell,n} = k)P(N_{\ell+1,n} = m | N_{\ell,n} = k)$. But in these cases, note that $P(N_{\ell+1,n} = m | N_{\ell,n} = k) = 0$ unless $k=m-1$. Hence we have \begin{eqnarray}
		P(N_{\ell+1,n} = m) &=& P(N_{\ell,n} = m-1)P(N_{\ell+1,n} = m | N_{\ell,n} = m-1)\nonumber\\ &=& P(N_{\ell,n} = m-1)p_{m-1}=\frac{\prod\limits_{j=1}^{m-2}p_j}{\sum\limits_{k=1}^\infty \prod\limits_{j=1}^{k-1} p_j}(p_{m-1}) = \frac{\prod\limits_{j=1}^{m-1}p_j}{\sum\limits_{k=1}^\infty \prod\limits_{j=1}^{k-1} p_j},\nonumber
	\end{eqnarray}
	which completes the induction step.

	To apply Lemma \ref{lem:bounded}, it suffices to show that \[
	T_n \le b\star N_{n,n}
	\]
	because $N_{n,n}$ all have the same distribution, and $T_n$ is a nondecreasing sequence. Then we conclude $P(T_n \text{ converges})>0$, and Kolmogorov's 0-1 law implies almost sure convergence. Define $T_{n}(\ell)$ by \[
	T_n(\ell) = \et_{k=n-\ell+1}^n A_{k}\,\,.
	\]
	Then $T_n(n) = \et_{k=1}^n A_{k} = T_n$. We can prove inductively (on $\ell$) that $T_n(\ell) \le b\star N_{\ell,n}$. This is trivial for $\ell=0$ because $T_n(0)=1$. Then observe:\[
	T_n(\ell+1) = \et_{k=n-\ell}^n A_{k} = (A_{n-\ell})^{\et_{k=n-\ell+1}^n A_{k}} = (A_{n-\ell})^{T_n(\ell)}\le (A_{n-\ell})^{b\star N_{\ell,n}}.
	\]
	If $A_{n-\ell} \le b^{\frac1{b\star N_{\ell,n}}}$ then by definition $N_{\ell+1,n} =1$ and thus $T_n(\ell+1) \le b = b\star N_{\ell+1,n}$. Otherwise, we have $N_{\ell+1,n} = N_{\ell,n} +1$, but still $A_{n-\ell}\le b$, and hence $T_{n}(\ell+1) \le b^{b\star N_{\ell,n}} = b\star N_{\ell+1,n}$. Thus we have for all $\ell$:\[
	T_n(\ell) \le b\star N_{\ell,n},
	\]
	and in particular $T_n=T_n(n)$ is a nondecreasing sequence with $b\star N_{n,n}$ as an upper bound whose distribution is constant in $n$, and therefore Lemma \ref{lem:bounded} implies that $T_n$ converges almost surely.\\
	\subsubsection{Convergence of $T_n$ implies summability of (\ref{eq:sum_check})}
	This time we assume $\ilim n T_n$ converges, and use its limiting distribution to construct a random variable $U_n$ on $\nn$ whose distribution's existence is tied to (\ref{eq:sum_check}) in a similar way to $N_{n,n}$ from the other direction. The random variable $U_n$ will also measure runs of large values of $A_i$, but in a different way. The construction hinges on the fact that $b\star n\rightarrow \infty$, which follows from the assumption that $b>e^{1/e}$. We also will find it convenient to extend the i.i.d. sequence $A_n$ to the non-positive integer indices. Using almost sure convergence of $T_n$ and the fact that $(A_n)_{n\in\zz}$ is i.i.d., we have\[
	S_n = \et_{k=0}^\infty A_{n-k}
	\]
	is almost surely finite, satisfies $S_{n+1}=(A_{n+1})^{S_n}$, and is identically distributed. We also observe that $S_n$'s distribution has unbounded support. Define $U_n = \min\{U\in \nn \mid (b\star U)^3> S_n\}$. This is well-defined because $b\star n$ is unbounded, also note that its distribution is constant in $n$. We ultimately will use the distribution of $U_1$ to put an finite upper bound on $\sum\limits_{k=1}^\infty\prod\limits_{j=1}^{k-1}p_j$. Choose $K\in \mathbb{N}$ such that $b\star(K-1) \ge 3$ and such that $P(U_1 = K+1)>0$. This exists because $\ilim k b\star k =\infty$, and the support of $U_1$'s distribution is unbounded. In order to tie $U_1$'s distribution to our desired sum, we will need to define an auxiliary sequence of random integers $(L_n)_{n\ge 0}$ recursively by $L_0 = 0$ and \[
	L_{n+1} = \begin{cases}\max\left\{L\in \{0,\dots,K\} \hspace{2pt}\vert\hspace{2pt} (b\star L)^3 \le S_{n+1} \right\} & L_n<K\\
		1+L_n & L_n \ge K \text{ and }A_{n+1}>b^{\frac1{b\star(L_n-1)}}\\
		0 & L_n \ge K \text{ and }A_{n+1}\le b^{\frac1{b\star(L_n-1)}}\end{cases}
	\]
	We claim $L_n\le U_n$ for all $n$, which can be proven by inducting $n$. The base case is trivial as $L_0 = 0 \le 1\le U_0$. It is similarly trivial that $L_{n+1}\le U_{n+1}$ when $L_n\ge K$ and $A_{n+1} \le b^{\frac1{b\star(L_n-1)}}$. If $L_n<K$, then \[
	(b\star L_{n+1})^3\le S_{n+1} < (b\star U_{n+1})^3
	\]
	from the definitions, which implies $L_{n+1}\le U_{n+1}$. In the remaining case, i.e. when $L_n \ge K$ and $A_{n+1}>b^{\frac1{b\star(L_n - 1)}}$, we claim that $U_{n+1}\ge U_n+1$. The inductive hypothesis $L_n\le U_n$ implies that in this case, $U_n\ge K$ and $A_{n+1} >  b^{\frac1{b\star(U_n-1)}}$, whence we can compute \begin{eqnarray}
		&& S_n\ge(b\star (U_n-1))^3 \nonumber\\
		\implies &&S_{n+1}=(A_{n+1})^{S_n}\ge(A_{n+1})^{(b\star (U_n-1))^3} \nonumber\\
		\implies &&S_{n+1}>b^{\frac{(b\star (U_n-1))^3}{b\star(U_n-1)}} \nonumber\\
		\implies &&S_{n+1}>\left(b^{b\star (U_n-1)}\right)^{b\star (U_n-1)} \nonumber\\
		\implies &&S_{n+1}>(b\star U_n)^{b\star (U_n-1)} \ge (b\star U_n)^{b\star (K-1)}\nonumber\\
		\implies &&S_{n+1}>(b\star U_n)^3 \nonumber
	\end{eqnarray} 
	This implies that $U_{n+1}\ge U_n+1$ by the definition of $U_n$. Therefore, if $U_n\ge L_n$ and $L_n\ge K$ and $A_{n+1} >  b^{\frac1{b\star(L_n-1)}}$, we can conclude $U_{n+1} \ge U_n+1$ and hence \[
	U_{n+1} \ge U_n+1\ge L_n + 1 = L_{n+1}
	\]
	in this case. Therefore we have by induction that $L_n\le U_n$ for all $n$.\\
	The distribution of $L_n$ will be somewhat similar to that of $N_{\ell,n}$ from the first part of the proof. Observe that for $m\in\nn$ we have $P(L_n = K+m\mid L_{n-1} = j) = 0$ except when $j=K+m-1$. Hence\begin{eqnarray}
		P(L_n = K+m) &=& \sum\limits_{j=0}^\infty P(L_{n-1} = j)P(L_n=K+m\mid L_{n-1} = j) \nonumber\\
		&=& P(L_{n-1} = K+m-1)P(L_n=K+m\mid L_{n-1} = K+m-1)\nonumber\\
		&=& P(L_{n-1} = K+m-1)P\left(A_n > b^{\frac1{b\star(K+m-1-1)}}\right)\nonumber\\
		&=& P(L_{n-1} = K+m-1) p_{K+m-2}.\nonumber
	\end{eqnarray}
	By repeatedly applying this fact, we obtain $ P(L_n = K+m) = P(L_{n-m}=K) \prod\limits_{j=K-1}^{K+m-2} p_j$ and therefore \[
	P(L_n\ge K) \ge P(L_n\in [K,K+n]) = \sum_{m=0}^n P(L_{n-m}=K) \prod_{j=K-1}^{K+m-2} p_j
	\]
	Observe that $L_n=K$ is equivalent to $(b\star K)^3\le S_n < (b\star (K+1))^3$, which in turn is equivalent to $U_n=K+1$. Therefore, using the fact that $U_n$ is constant in distribution and $U_n\ge L_n$,\begin{eqnarray}
		P(U_1\ge K) &=& P(U_n\ge K) \ge P(L_n\ge K) \ge\sum_{m=0}^n P(L_{n-m}=K) \prod_{j=K-1}^{K+m-2} p_j \nonumber\\&=& \sum_{m=0}^n P(U_{n-m}=K+1) \prod_{j=K-1}^{K+m-2} p_j = P(U_1=K+1) \sum_{m=0}^n  \prod_{j=K-1}^{K+m-2} p_j\nonumber
	\end{eqnarray}
	and thus \[
	\sum_{m=0}^n  \prod_{j=K-1}^{K+m-2} p_j\le \frac{P(U_1\ge K)}{P(U_1 = K+1)} < \infty 
	\]
	for all $n$ and hence \[
	\sum_{m=0}^\infty  \prod_{j=K-1}^{K+m-2} p_j \le \frac{P(U_1\ge K)}{P(U_1 = K+1)} < \infty 
	\]
	Convergence of this series is equivalent to convergence of the series from (\ref{eq:sum_check}).

	\subsubsection{Equivalence of finiteness of (\ref{eq:sum_check}) and $\ee\left[\inf\{n: A_n \le e^{\frac 1 { e\star n}}\}\right] $ }	
	We start by observing that for any $c>1$ \[\ee \left[\inf\{n: A_n \le c^{\frac 1 { c\star n}}\}\right] = \sum_{k=1}^\infty \prod_{j=1}^{k-1}  P(A_1>c^{\frac1{c\star j}})\]
    For $c=b$ this is simply the summation in equation (\ref{eq:sum_check}) . Therefore it suffices to prove that $\ee\left[ \inf\{n: A_n \le b^{\frac 1 { b\star n}}\}\right]$ is infinite if and only if $\ee \left[\inf\{n: A_n \le e^{\frac 1 { e\star n}}\}\right]$ is. Of course, this would follow if, for some fixed $k\in \mathbb N$ depending only on $b$, \begin{equation}\label{eq:inf}
	P(A_1>b^{\frac1{b\star ({j+k})}})\le P(A_1>e^{\frac1{e\star j}}) \le P(A_1>b^{\frac1{b\star ({j-k})}}) 
    \end{equation}
    for all $j\in \nn$.
	To show this, we apply this lemma on the growth rate of power towers:
	\begin{lemma}\label{lem:growth}
		Suppose $s_n$, $t_n$ are two sequences on $(e^{\frac1e},\infty)$ and $u,c\in(e^{1/e},\infty)$ such that \begin{eqnarray*}
			s_n \le u\star n &\hspace{6pt}\mbox{and}\hspace{6pt}& t_n \le u\star n\\
			s_n \ge c &\hspace{6pt}\mbox{and}\hspace{6pt}& t_n \ge c.
		\end{eqnarray*}
		Let $T_n = t_n^{t_{n-1}^{{\iddots^{t_1}}}}$ and $S_n = s_n^{\iddots^{s_1}}$. Then, for each $\lambda\in(0,\infty)$ there exists a constant $k\in \nn$ such that \[
		S_{n-k} \le \lambda T_n \le S_{n+k}
		\]
		for all $n >k$.
	\end{lemma}
	
	When $b>e^{1/e}$, we can take $s_n\equiv b$, $t_n\equiv e$, and $\lambda=\frac1{ln b}$ and conclude that there exists $k\in\nn$ such that \[
	e^{\frac1{e\star(n+k)}}\le b^{\frac1{b\star n}} \le e^{\frac1{e\star(n-k)}},
	\]
	from which (\ref{eq:inf}) follows immediately. This completes the proof of Theorem \ref{th:main} case 2.\qed\\
	\emph{Proof of Lemma \ref{lem:growth}.} We first show the result for $t_n = u$ and $T_n = u\star n$: Pick $\tau > 1$. Observe that \[
	\ilim j u\star j = \ilim j \frac{(u\star j)^{1-\frac1\tau}}{u\star(j-1)} = \infty.
	\]
	Therefore we can choose $k_\tau\in\nn$ large enough such that \[
	u\star(j+1) > s_1^\tau\hspace{7pt}\text{and}\hspace{7pt}
	\frac{(u\star j)^{1-\frac1\tau}}{u\star(j-1)} > \tau
	\]
	for all $j\ge k_\tau$.\\
	We will show inductively that \begin{equation}\label{eq:base_case}
		u\star n > S_{n-k_\tau}^\tau
	\end{equation}
	for all $n>k_\tau$. We already have assumed the base case, i.e. $n=k_\tau+1$. Suppose that (\ref{eq:base_case}) holds for some $n>k_\tau$. Then we make a lengthy computation:\begin{eqnarray*}
		\frac{(u\star n)^{1-\frac1\tau}}{u\star(n-1)} &>& \tau\\
		\frac{(u\star n)^{1-\frac1\tau}}{u\star(n-k_\tau)} &>& \tau\\
		\frac{u\star n}{u\star (n-k_\tau)} &>& \tau \left(u\star n\right)^{\frac1\tau}\\
		u\star n &>& \tau \left(u\star n\right)^{\frac1\tau}(u \star(n-k_\tau))\\
		u\star n &>& \tau({u\star (n-k_\tau)}) S_{n-k_\tau}\\
		u^{u\star n}&>& u^{\tau({u\star (n-k_\tau)}) S_{n-k_\tau}} = \left(\left(u^{u\star(n-k_\tau)}\right)^{S_{n-k_\tau}}\right)^\tau\\
		u\star(n+1)&>& \left(\left(u\star(n+1-k_\tau)\right)^{S_{n-k_\tau}}\right)^{\tau}\\
		u\star(n+1)&>&\left(\left(s_{n+1-k_\tau}\right)^{S_{n-k_\tau}}\right)^{\tau} = (S_{n+1-k_\tau})^\tau.
	\end{eqnarray*}
	Thus, by induction, we have $u\star n > S_{n-k_\tau}^\tau$ for all $n>k_\tau$.\\
	For the other inequality, pick $\beta<\min\{\frac{\ln c}{\ln u},1\}$, then chose $k_\beta\in\nn$ large enough so that \[
	1 = u\star 0 < \beta S_{k_\beta}
	\]
	and \[
	\frac{-\ln\beta}{\ln c - (\ln u) \beta} < S_{k_\beta}.
	\]
	Then we show inductively that \begin{equation}\label{eq:base_case2}
		u\star n < \beta S_{n+k_\beta}
	\end{equation}
	for all $n$. We have the base case ($n=0$) already. Note that the second inequality implies \[
	(\ln u)\beta S_{n+k_\beta} < \ln \beta + (\ln c) S_{n+k_\beta}
	\]
	for all $n\ge 0$. Assuming $u\star n < \beta S_{n+k_\beta}$, this implies \begin{eqnarray*}
		(\ln u)(u\star n) &<& \ln \beta + (\ln c) S_{n+k_\beta}\\
		u\star (n+1) &<& \beta c^{S_{n+k_\beta}}\\
		u\star (n+1) &<& \beta (s_{n+k_\beta+1})^{S_{n+k_\beta}}\\
		u\star (n+1) &<& \beta S_{n+k_\beta+1},
	\end{eqnarray*}
	thus we have that (\ref{eq:base_case2}) holds for all $n\ge 0$.\\
	It is apparent that $\ilim n \frac{u\star(n+1)}{u\star n}=\infty$, hence there exists $k_\lambda$ such that \[
	u\star(n-1)\le \lambda (u\star n) \le u\star (n+1)
	\]
	for all $n\ge k_\lambda$. Therefore, if we take $k = 1+\max\{k_\lambda,k_\tau,k_\beta\}$, we will have, for any $n\ge k$ that \begin{eqnarray}
	S_{n-k}&\le& S_{n-1-k_\tau}<S_{n-1-k_\tau}^\tau< u\star(n-1)\nonumber\\
	&\le& \lambda(u\star n)\le u\star(n+1)<\beta S_{n+1+k_\tau}< S_{n+1+k_\tau}\le S_{n+k},\nonumber
	\end{eqnarray}
	as desired.\\
	Now to prove the general case. Suppose $s_n$ and $t_n$ are two sequences satisfying the theorem's bounds. As we have shown, there exist $k_1,k_2\in \nn$ such that \begin{eqnarray*}
		T_{n-k_1}<u\star n < T_{n+k_1} &\hspace{6pt}\text{and}\hspace{6pt}& S_{n-k_2}<u\star n < S_{n+k_2}
	\end{eqnarray*}
	for all $n>\max(k_1,k_2)$. Thus we have \[
	S_{n-k_2-k_1} < u\star(n-k_1) < T_n < u\star(n-k_2) < S_{n+k_1+k_2},
	\]
	so choosing $k=k_1+k_2$ we obtain the desired result.\qed

\subsection{Cases 3 and 4 ($a<1$)} Because $a<1$, $T_n$ is not a monotonic sequence. A different approach is needed. Our proofs hinge on the following theorem.
\begin{theorem}\label{th:monotonic_compositions}
	Let $f_i:[a,b] \rightarrow [a,b]$ be a sequence of i.i.d. randomly distributed non-decreasing functions, where $[a,b]$ is a closed interval in the extended real line $\rr\cup\{\pm\infty\}$. Furthermore, assume that $P(f_0(a) \ge f_1(b))>0$. We define $F_n:[a,b]\rightarrow[a,b]$ by \[
	F_n(x) = (f_0 \circ f_1 \circ \cdots \circ f_n)(x).
	\]
	Then $F_n(x)$ almost surely converges, and its limit is independent of $x$.
\end{theorem}
Unlike the bulk of previous work on iterated random functions, this theorem makes no continuity assumptions. Furthermore, unlike Lemma \ref{lem:bounded} and the ``contracting on average" results in Diaconis and Freedman, this theorem has no obvious non-probabilistic analogue. It does bring to mind the fact that any non-decreasing function mapping a closed interval to itself must have a fixed point (a consequence of Tarski's lattice-theoretical fixed point theorem \cite{article:3}), but Tarski's theorem has no implications about the dynamics of such a function. Indeed, the iterates of a non-decreasing function of a closed interval need not converge to a fixed point in general.

\begin{proof}Without loss of generality, we suppose $[a,b]$ is finite (in the other case, we can use a bounded increasing function to convert $[a,b]$ to a finite interval). Observe that $F_n(a)$ must form a non-decreasing sequence because $f_n(a)\ge a$, and similarly, $F_n(b)$ forms a non-increasing sequence, hence both of these converge almost surely. Call $F_\infty(a)$ and $F_\infty(b)$ their respective limits. Let $f_\omega$ and $f_{\omega+1}$ be i.i.d. with the same distribution as the $f_i$'s but independent of that sequence. We let $X_n = F_n(f_\omega(a))$ and $Y_n = F_n(f_{\omega+1}(b))$. Note that $X_n \stackrel{d}{=} F_{n+1}(a)$ and $Y_n \stackrel{d}{=} F_{n+1}(b)$. We thus have $X_n\toindist F_\infty(a)$ and $Y_n\toindist F_\infty(b)$. We furthermore have convergence in probability: $X_n \ge F_n(a)$, hence $X_n -F_n(a) \ge 0$, and by boundedness of $[a,b]$, we also have $\ee(X_n - F_\infty(a))\rightarrow 0$, which implies $X_n\toinp F_\infty(a)$. Similarly $Y_n\toinp F_\infty(b)$. \\
	Since $F_n$ is always nondecreasing, we have \begin{equation*}
		P(Y_n \le X_n) = P{\big[}F_n(f_{\omega+1}(b)) \le F_n(f_{\omega}(a)){\big]} \ge P\left(f_{\omega+1}(b) \le f_{\omega}(a)\right) > 0.
	\end{equation*}
	Hence we have a constant, nonzero lower bound on $P(Y_n\le X_n)$. By convergence in probability of $X_n$ and $Y_n$ to $F_\infty(a)$ and $F_\infty(b)$, respectively, we can conclude the same inequality for their limits \[P(F_\infty(b) \le F_\infty(a)) \ge P(f_{\omega+1}(b) \le f_{\omega}(a)) > 0.\]
	It is clear from the definition that $F_\infty(b) \ge F_\infty(a)$ almost surely, so this implies that there is a nonzero probability $F_\infty(b) = F_\infty(a)$. By monotonicity of $F_n$ for each $n$, $F_\infty(b) = F_\infty(a)$ is equivalent to $F_n$ converging uniformly to a constant function. Thus, we have a nonzero probability of $F_n$ converging to a constant function. Since convergence of $F_n$ to a constant function is a tail event, this implies that the convergence is almost sure since its probability is nonzero.\end{proof}

We will apply Theorem \ref{th:monotonic_compositions} to both cases 3 and 4 of Theorem \ref{th:main}.
\begin{theoremcase}
	Let $B:[0,1]\rightarrow[0,1]$ be defined by \[
	B(x) =\begin{cases} x& x\ge e^{-e}\\\exp\left(W\left(\frac1{\ln x}\right)\exp\left(-\frac1{W\left(\frac1{\ln x}\right)}\right)\right) & x<e^{-e}\end{cases}.
	\]
	If $b\le1$ and $a<B(b)$, then $T_n$ converges almost surely. Conversely, if $b<e^{-e}$ and $a\ge B(b)$, then $T_n$ diverges by oscillation, with the even an odd subsequences converging to distinct limits.
\end{theoremcase}
\begin{proof}
	We will show that, when $a<B(b)$, the functions $x\rightarrow \et_{i=2kn+1}^{2 k(n+1)}A_i ; x$ satisfy the conditions of Theorem \ref{th:monotonic_compositions} for some $k$. When $a\ge B(b)$, we will show that there is almost surely a non-zero lower bound on the difference between $T_{2n}$ and $T_{2n+1}$. 
	
	Before diving into the proof, we note that for $c\in[0,1]$ function $x\to c^x$ is a non-increasing function from $[0,1]$ to itself, so we can restrict all our attention to this interval as $T_n$ will always be in that interval, and furthermore either $T_n$ converges in the limit, or it diverges by oscillation. We also have that $T_{2n}$ is a non-increasing sequence, and $T_{2n+1}$ is a non-decreasing sequence, which follows from the fact that composing two exponential maps with base less than 1 forms a non-decreasing function. We furthermore observe that $\et_{k=1}^n c_k$ is nondecreasing in its even terms and nonincreasing in its odd terms.
	
	Now, define the alternating power tower function $\at n(x,y)$ recursively by $\at 0(x,y)=1$ and \[
	\at n(x,y) = x^{\at {n-1}(y,x)},
	\]
	so called because it makes a power tower that alternates between $x$ and $y$, e.g.\[
	\at 5(x,y) = x^{y^{x^{y^x}}}.
	\]
	These have the monotonicity properties that for each $n$, $\at n (x,y)$ is increasing in $x$ and decreasing in $y$. Furthermore, for fixed $x,y$, we have that $\at{2n-1}(x,y)$ is an increasing sequence and $\at{2n}(x,y)$ is a decreasing sequence, so we can also define the functions\begin{eqnarray}
		\ate(x,y) &=& \ilim n \at{2n}(x,y)\nonumber\\
		\ato(x,y) &=& \ilim n \at{2n-1}(x,y)\nonumber
	\end{eqnarray}
	which obey the relations \[
	\ate(x,y)=x^{\ato(y,x)} \hspace{24pt} \ato(x,y)=x^{\ate(y,x)}
	\]
	and furthermore, both are fixed points of the map $t\to x^{y^t}$. In fact, $\ate(x,y)$ is the largest such fixed point, and $\ato(x,y)$ is the smallest. Baker and Rippon \cite{ARTICLE:7} have shown that there exists a function $B(x)$ such that, for $y\le x$,  $\ate(x,y)=\ato(x,y)$ if and only if $0<y\le B(x)$. Here, we have found its closed form in terms of the Lambert W-function.
	
	This case of Theorem \ref{th:main} follows immediately from these two lemmas:
	
	\begin{lemma}\label{lem:case_3_part_1}
		If $0\le a<b\le 1$ and there exists $k\in\nn$ such that $\at{2k-1}(b,a)>\at{2k}(a,b)$, then $T_n$ converges almost surely. Conversely, if $\at{2k-1}(b,a)\le \at{2k}(a,b)$ for all $k$, then $T_n$ almost surely diverges by oscillation.
	\end{lemma}
	\begin{lemma}\label{lem:B}
		If $0\le a < B(b)$, then $\at{2k-1}(b,a)>\at{2k}(a,b)$ for some $k\in\nn$. Conversely, if $B(b)\le a\le b$, then $\at{2k-1}(b,a)\le \at{2k}(a,b)$ for all $k\in\nn$.
	\end{lemma}
	\emph{Proof of Lemma \ref{lem:case_3_part_1}.} 
	Suppose $\at{2k-1}(b,a)>\at{2k}(a,b)$. Then define $f_n:[0,1]\rightarrow[0,1]$ by \[f_n(x) = \et_{i={2nk+1}}^{2(n+1)k} A_i ; x.\]
	We note that $(f_n)_{n\in\nn}$ is a random i.i.d. sequence of increasing functions on a closed interval, so we will attempt to apply Theorem \ref{th:monotonic_compositions}. By continuity, we have that for some $\epsilon>0$ sufficiently small, $\at{2k-1}(b-\epsilon,a+\epsilon)>\at{2k}(a+\epsilon,b-\epsilon)$. Also by assumption, we have $P(A_i \in[a,a+\epsilon])$ and $P(A_i\in[b-\epsilon,b])$ are positive for all $i$, hence \begin{alignat}{3}
	P(f_0(0)\ge f_1(1)) &=P\left(\et_{i={1}}^{2k} A_i ; 0 \ge \et_{i={2k+1}}^{4k} A_i ; 1\right) = P\left(\et_{i=1}^{2k-1}A_i \ge \et_{i={2k+1}}^{4k} A_i\right)\nonumber\\
	&\ge P\left(\et_{i=1}^{2k-1}A_i > \at{2k-1}(b-\epsilon,a+\epsilon)\, \mathrm{ and }\et_{i={2k+1}}^{4k} A_i <\at{2k}(a+\epsilon,b-\epsilon))\right)\nonumber\\ &> 0.\nonumber
	\end{alignat}	
	Thus we can apply the result of Theorem \ref{th:monotonic_compositions} to obtain \begin{eqnarray}
	\ilim n f_1\circ f_2\circ  \cdots \circ f_n(0) &=& \ilim n f_1\circ f_2\circ  \cdots \circ f_n(1)\nonumber\\
	 \ilim n \et_{i=1}^{2nk}A_i;0 &=& \ilim n \et_{i=1}^{2nk}A_i;1\nonumber\\
	 \ilim n T_{2kn-1} &=& \ilim n T_{2kn}\nonumber.
	\end{eqnarray}
	To complete the proof, we use the fact that both the even- and odd-indexed power towers converge. We conclude \[
	\ilim n T_{2n-1} = \ilim n T_{2kn-1} = \ilim n T_{2kn} = \ilim n T_{2n},
	\]
	hence $T_n$ converges, as desired. 
	
	Conversely, we suppose $\at{2k-1}(b,a)\le \at{2k}(a,b)$ for all $k$, which implies $\ato(b,a) \le \ate(a,b)$. Let $c\in(a,b]$ such that $P(A_1\ge c)>0$. Therefore, we have the following with nonzero probability: \begin{eqnarray}
	T_{2k}-T_{2k-1} &=& (A_1)^{\et_{i=2}^{2k}A_i}-\et_{i=1}^{2k-1}A_i\ge (A_1)^{\at{2k-1}(b,a)}-\at{2k-1}(b,a)\nonumber\\
	&\ge& (A_1)^{\ato(b,a)}-\ato(b,a)\ge c^{\ato(b,a)}-\ato(b,a)\nonumber\\ &>& a^{\ato(b,a)}-\ato(b,a) = \ate(a,b)-\ato(b,a) \ge 0\nonumber
	\end{eqnarray}
	hence we have a nonzero probability that $\ilim n T_n$ does not converge, since there is a nonzero probability that the limit of $T_{2k}$ is strictly greater than the limit of $T_{2k-1}$. Since convergence is a tail event, Kolmogorov's 0-1 law implies that $\ilim n T_n$ converges with probability $0$, proving the lemma.\\
	\emph{Proof of Lemma \ref{lem:B}.}  
	Baker and Rippon \cite{ARTICLE:7} have proven that $\ilim n \at n (a,b)$ converges if and only if $\phi_{a,b}(x)= a^{b^x}$ has exactly one fixed point $c$ such that $|\phi_{a,b}'(c)|\le 1$. Furthermore, they showed for each $b$, there exists a constant $a_1$ such that this happens for all $a\le a_1$ and for no $a\in(a_1,b)$. Furthermore, $\phi_{a_1,b}(x)$ has a fixed point $c$ such that $\phi'_{a_1,b}(c) = 1$. 
 
    If $\at n(a,b)$ converges, by monotonicity and the fact that $a\ne b$, this implies $$
	\ilim n \at n(b,a) > \ilim n \at n (a,b),
	$$
	and in particular that $\at{2n-1}(b,a) > \at {2n}(a,b)$ for some $n$. Conversely, if $\ilim n \at n(a,b)$ does not converge, we must have $a\in(a_1,b)$, and by continuity of $a^{b^x}$ in all variables and its monotonicity properties, we must have that $\ate(a,b)$ is larger than the fixed point of $x\rightarrow b^x$ for all $a\in (a_1,b)$. This implies that $\ate(a,b)^{1/{\ate(a,b)}} \ge b$, and thus $\ate(a,b) \ge b^{\ate(a,b)} = \ato(b,a)$. Hence, for all $n$, we must have $\at{2n}(a,b)\ge \at{2n-1}(b,a)$.
 
	Therefore it suffices to show that $\phi_{a,b}(c) = c$ and $|\phi_{a,b}'(c)| \le 1$ has exactly one solution if and only if $0\le a< B(b)$. Baker and Rippon also proved that for each $b\ge e^{-e}$, this holds for all $a\in[0,b]$, and for each $b<{e^{-e}}$, there exists a constant $a_1\in(0,b)$ such that this holds for $a<a_1$ and fails for $a\in [a_1,b]$. Furthermore, they proved that there is a fixed point $c$ of $\phi_{a_1,b}$ such that $\phi'_{a_1,b}(c)=(\log a_1)(\log b) b^c a_1^{b^c} = 1$. Therefore, it suffices to prove that $a=B(b)$ is the only number in $(0,b)$ that solves the system \begin{eqnarray}
	a^{b^c} &=&c\nonumber\\
	(\log a) (\log b) b^c a^{b^c} &=& 1\nonumber.
	\end{eqnarray}
	Solving this system is not obviously possible in any kind of closed form. Using the Lambert $W$-function, it becomes straightforward. Start by rearranging the second equation to \begin{eqnarray}
	\log(a^{b^c}) a^{b^c} &=& \frac1{\log b}\nonumber\\
	\log(c) c &=& \frac1{\log b}\nonumber\\
	\log c &=& W\left(\frac1{\log b}\right)\nonumber
	\end{eqnarray}
	where this can be either of the two real branches of the $W$ function. Using the fact that $\exp(W(x)) = \frac x {W(x)}$ for any branch of the $W$ function, this gives us $c = \frac1{(\log b) W\left(\frac1{\log b}\right)}$, and thus $b^c = \exp\left(\frac1{W\left(\frac1{\log b}\right)}\right)$. Thus, the requirement $a=c^{b^{-c}}$ gives us the desired form for $a=B(b)$, and the requirement that $a\le b$ implies that we must take the principle branch of the $W$-function. 
	\end{proof}

We can do a similar technique, using Theorem \ref{th:monotonic_compositions}, to solve the case of $a<1<b$.
\begin{theoremcase}
	If $a<1<b\le\infty$, then $T_n$ converges almost surely.
\end{theoremcase}
\begin{proof}
	We obviously have $P(A_n<1)>0$, therefore we have a.s. that $A_n<1$ infinitely often. Let $(n_i)_{i\ge 0}$ be the set of indices where $A_n<1$ (note that $n_i$ is an increasing sequence of random variables, and also $n_{i+1}  - n_i$ is i.i.d.). Define $$
	f_i(t) = \et_{k=n_{2i}}^{n_{2i+2}-1} A_k ; (t).
	$$
	By letting $f_i(\infty)=\ilim t f_i(t)$, we can consider this as a function from $[0,\infty]$ to $[0,\infty]$ and apply Theorem \ref{th:monotonic_compositions}. (The limit is well defined because the limit of a monotonic function on $[0,\infty)$ either converges or diverges to infinity.) To show that this theorem is applicable, we must show that $f_i(t)$ are increasing in $t$, i.i.d., and $P(f_0(0)\ge f_1(\infty))>0$. 
	
	Firstly, $f_i$ is increasing because $\{A_{n_{2i}},\dots, A_{n_{2i+2}-1}\}$ contains exactly two elements less than one. Secondly, $f_i$ are i.i.d. because the three sequences $(A_{n_i})_{i\ge 1}$, $\{A_n\mid A_n\ge 1\}$, and $(n_i)_{i\ge 1}$ are independent and each individually is i.i.d. Thirdly, we must show $P(f_0(0)\ge f_1(\infty))>0$. Observe that \[f_0(0) = \et_{k=n_{2i}}^{n_{2i+2}-1} A_k ; 0 = \et_{k=n_{2i}}^{n_{2i+2}-2} A_k	\]
	and \[
	f_1(\infty) = \ilim t f_1(t) = \ilim t  \et_{k=n_{2i}}^{n_{2i+2}-1} A_n ; (t) = \ilim x \et_{k=n_{2i}}^{n_{2i+1}} A_k ; x = \et_{k=n_{2i}}^{n_{2i+1}-2} A_k\,.
	\]
	Next, choose $c\in[a,1)$ such that both $P(A_n\in [c,1))$ and $P(A_n\in [a,c])$ are nonnegative. This exists by the definition of $a$. Thus the event that $A_1\in[c,1)$ and $A_2,A_3\in [a,c]$ and $A_4>1$ has nonzero probability. In this event, we have $n_i=i+1$ for $i\le 2$ and $n_3 >4$, and the following inequalities:\begin{eqnarray*}	f_0(0) &=& \et_{k=n_0}^{n_2-1} A_k = \et_{k=1}^{2} A_k = (A_1)^{A_2} \ge c^{A_2} \ge c^c
    \\
	f_1(\infty) &=& \et_{k=n_{2}}^{n_{3}-1} A_k = (A_3)^{\et_{k=4}^{n_3-1} A_k} \le A_3 \le c.
	\end{eqnarray*}
	Because $c^c \ge c$, we therefore conclude that \begin{eqnarray*}
	&&P\left(f_0(0) \ge f_1(\infty)\right) \ge P\left(f_0(0) \ge c^c\, \mathrm{ and } \,f_1(\infty) \le c\right)\\&&\ge P\left(A_1\in[c,1) \text{ and }A_2,A_3\in[a,c]\text{ and }A_4>1\right)>0.\end{eqnarray*}
	Therefore we can apply Theorem \ref{th:monotonic_compositions} to \[
	F_i(x) = f_0\circ f_1 \circ \cdots  \circ f_i(x) = \et_{k=n_0}^{n_{2i+2}-1}A_k;x 
	\]
	and find that $\ilim i F_i(x)$ converges almost surely to a limit that does not depend on $x$. In particular, $\ilim i F_i(0) = \ilim i F_i(\infty)$. From the definition, we note that $f_i(t)$ is always a number less than $1$ to a nonnegative power, hence it is bounded above by $1$, and therefore $\ilim i F_i(x)$ is finite. To complete the proof, we note what this means for $T_n$. Let $i(n) = \max\{i\mid n_{2i+2}-1\le n\}$. Then \[
	T_n = \et_{k=1}^n A_k = \et_{k=1}^{n_0-1}A_k;\et_{k=n_0}^{n_{2i(n)+2}-1}A_k ;\et_{k=n_{2i(n)+2}}^{n}A_k =\et_{k=1}^{n_0-1}A_k;\left(F_{i(n)}\left(\et_{k=n_{2i(n)+2}}^{n}A_k \right)\right),
	\]
	which implies \[
	\et_{k=1}^{n_0-1}A_k;\left(F_{i(n)}\left(0 \right)\right)\le T_n \le \et_{k=1}^{n_0-1}A_k;\left(F_{i(n)}\left(\infty \right)\right).
	\]
	The upper and lower bounds almost surely converge to the same finite value, and hence $T_n$ converges almost surely.
\end{proof}

\subsection{Proof of Corollaries}
The corollaries follow trivially from Theorem \ref{th:main} except in the cases where $b>e^{\frac1e}$ and $a\ge 1$, where we have the claim that $T_n$ diverges if $a>1$, if $\ee \left[\log^\star(\frac1{A_1-1})\right] < \infty$, or if $\liminf\limits_{t\rightarrow 1^+}\left[P(A_1\le t)\log^\star\left(\frac1{t-1}\right)\right]>1$, and $T_n$ converges if $\limsup\limits_{t\rightarrow 1^+}\left[P(A_1\le t)\log^\star\left(\frac1{t-1}\right)\right]<1$. Note that $a>1$ implies $\ee \left[\log^\star(\frac1{A_1-1})\right] < \log^\star(\frac1{a-1}) < \infty$. Also, observe that for any non-increasing function $f:[1,\infty)\rightarrow [0,\infty)$ and any random variable $X$ supported on a subset of $[1,\infty)$ with $\ee f(X) < \infty$, we have \begin{eqnarray*}
    0&\le& \limsup\limits_{t\rightarrow 1^+}\left[f(t)P(X\le t)\right] =  \limsup\limits_{t\rightarrow 1^+}\left[f(t)\int_1^\infty\mathbb{1}_{[1,t]}(x) dP_X(x) \right]=0,
\end{eqnarray*}
 using the dominated convergence theorem in the last step, which is applicable because $f(t)\mathbb{1}_{[1,t]}(x)\le f(x)$ for all $x$, and $f(x)$ is integrable with respect to $dP_X$. Therefore $\ee\left[\log^\star(\frac1{A_1-1})\right] <\infty$ implies $\limsup\limits_{t\rightarrow 1^+}\left[P(A_1\le t)\log^\star\left(\frac1{t-1}\right)\right]=0<1$, so we conclude that it suffices to prove\\
\begin{lemma}
    \[
	\liminf\limits_{t\rightarrow 1^+}\left[P(A_1\le t)\log^\star\left(\frac1{t-1}\right)\right]>1\implies \ee\left[\inf\{n : A_n \le e^{\frac1{e\star n}}\}\right] < \infty
	\] and
	\[\limsup\limits_{t\rightarrow 1^+}\left[P(A_1\le t)\log^\star\left(\frac1{t-1}\right)\right]<1 \implies \ee\left[\inf\{n : A_n \le e^{\frac1{e\star n}}\}\right] = \infty\].
\end{lemma}
\begin{proof}
Note that \[
    \ee\left[\inf\{n : A_n \le e^{\frac1{e\star n}}\}\right] =
        \sum_{n=1}^\infty \prod_{j=1}^{n-1} P(A_1 > e^{\frac1{e\star j}}) = \sum_{n=1}^\infty \prod_{j=1}^{n-1} \left(1-P(A_1 \le e^{\frac1{e\star j}})\right).
    \] 
    In this form, it is clear that having information on the asymptotic behavior of $P(A_1\le t)$ as $t$ goes to $1$ from above should give us information about the convergence of this sum. 
    
    First, consider the case of $P(A_1=1)>0$, which implies $\lim\limits_{t\rightarrow 1^+}\left[P(A_1\le t)\log^\star\left(\frac1{t-1}\right)\right]=\infty>1$ because $\lim\limits_{t\rightarrow 1^+}\log^\star\left(\frac1{t-1}\right)=\infty$, and also $P(A_1=1)>0$ implies $\ee(\inf\{n : A_n \le e^{\frac1{e\star n}}\}) \le \ee(\inf\{n : A_n =1\})< \infty$. Hence, we may proceed assuming that $P(A_1=1)=0$.
    
    Before proving either claim under this assumption, we want to replace $\log^\star(\frac1{A_1 -1})$ with a new variable $X=\log^\star(\frac1{\log A_1})$ as this will be easier to work with. To do this, we start by noting that $\abs{\frac{1}{\log x}-\frac1{x-1}} < 1$ on $(1,\infty)$, which may be shown by elementary calculus. We also know that for integers $n\ne m$, $\abs{e\star n -e \star m} > 1$, which implies that if $\abs{x-y} <1$ then $\abs{\log^\star(x)-\log^\star(y)} \le 1$. Therefore \[
	\liminf\limits_{t\rightarrow 1^+}\left[P(A_1\le t)\log^\star\left(\frac1{t-1}\right)\right] = \liminf\limits_{t\rightarrow 1^+}\left[P(A_1\le t)\log^\star\left(\frac1{\log t}\right)\right],
	\]
    and similarly the corresponding limsups are equal.
    
    Suppose that $\liminf_{t\rightarrow 1^+}\left[P(A_1\le t)\log^\star\left(\frac1{\log t}\right)\right]>1$. Then, there exist constants $c>1$ and $N\in\nn$ such that $j\ge N$ implies \[
    P(A_1\le e^{\frac1{e\star j}}) \ge \frac c{\log^\star\left(1/\log(e^{\frac1{e\star j}})\right) } = \frac c j.
    \]
    Thus, for $n>N$, we have \begin{alignat}3
        \prod_{j=N}^{n-1} \left(1-P(A_1 \le e^{\frac1{e\star j}})\right)& \le \prod_{j=N}^{n-1} \left(1-\frac c j\right) =\exp\left(\sum_{j=N}^{n-1} \log\left(1-\frac c j\right)\right)\le \exp\left(-\sum_{j=N}^{n-1} \frac c j\right)\nonumber\\
        &\le \exp\left(-c\log n + c\log N\right) = \frac{\exp(c\log N)}{n^c}\nonumber
    \end{alignat}
    Since $\sum_n \frac1{n^c} < \infty$, we conclude that $
    \sum_{n=1}^\infty \prod_{j=1}^{n-1} \left(1-P(A_1 \le e^{\frac1{e\star j}})\right) < \infty,
    $ as desired.
    Similarly, if $\limsup_{t\rightarrow 1^+}\left[P(A_1\le t)\log^\star\left(\frac1{\log t}\right)\right]<1$, the summands go to $0$ slower than $\frac1{n^{c}}$ for some $c\in(0,1)$, hence the series diverges in that case.
\end{proof}
\section{The Inverse Question}
It is natural to ask, which distributions can be represented as power towers of an i.i.d. sequence? More precisely, given a random variable $T$, does there exist an i.i.d. sequence $\{A_i\}_{i\in\nn}$ such that $T \eqindist (A_1)^{(A_2)^{(A_3)^{\iddots}}}
$? If the answer is yes, we will say $T$ has a tower distribution, and the distribution of $A_i$ we will call its inverse tower distribution. We will not answer this question in full generality, but Theorem \ref{th:inverse_problem} gives the answer in the case that $T = U^r$ for a uniform distribution $U$ and fixed $r\in\rr\setminus\{0\}$.
\begin{theorem}\label{th:inverse_problem}
	Let $U\sim\mathrm{Unif}[\alpha,\beta]$ for some $0\le \alpha<\beta$. Then $U^r$ has a tower distribution if and only if $1\in[\alpha,\beta]$ and $r\in[0,\frac1{1+\log \beta}]$. When these conditions hold and $r\ne 0$, the inverse tower distribution has distribution function given by \[
	F(x) = \begin{cases}\frac{x^{\frac{\beta^r}r}}{\beta}(1-\beta^r\log x) - \frac{\alpha}{\beta} + \frac \alpha\beta F\left(x^{\frac{\beta^r}{\alpha^r}}\right)&x\in[\alpha^{r/\beta^r} , \beta^{r/\beta^r})\setminus\{1\}\\
		0 & x < \alpha^{r/\beta^r}\\
		1 & x \ge \beta^{r/\beta^r}\\
		\frac{1-\alpha}{\beta-\alpha} & x = 1\end{cases}
	\]
    where we interpret $ \frac \alpha\beta F(x^{\frac{\beta^r}{\alpha^r}})$ as $0$ when $\alpha=0$.
\end{theorem}
This result suggests that the inverse question may be quite difficult for an arbitrary $T$, as it is not clear exactly what ``goes wrong" when the conditions are not met. In the case of $\alpha=0$ and $\beta=1$, this $F$ has a relatively simple form, and we get an interesting generalization of Example \ref{ex1}:
\begin{example}
	If $U\sim \mathrm{Unif}(0,1)$ then \[
	T= U^r
	\]
	has a tower distribution if and only if $r\in[0,1]$. When $r\in(0,1)$: Let $V_1,V_2,V_3$ be independent and uniform on $(0,1)$. Then \[
	A = \max(V_1^{\frac r {1-r}}, V_2V_3)
	\] 
	has the inverse tower distribution of $T$.
\end{example}
In our proof of Theorem \ref{th:inverse_problem}, we will use the following lemmas, each of which may be useful on its own for checking whether a distribution is a tower distribution: 
\begin{lemma}\label{lem:Letac}
	Suppose $X,Y$ are independent with almost sure bounds $X\in[0,e]$ and $Y\in[0,e^{\frac1e}]$. If $X \eqindist Y^X$, then $Y$ has a tower distribution, and its inverse tower distribution is $X$.
\end{lemma}
\begin{proof}
Let $X_i$ be an i.i.d. sequence with terms having the same distribution as $X$, but independent of it. By the second condition of Corollary \ref{cor:conv}, the power tower formed $T=\et_{i=1}^\infty X_i$ converges almost surely. This allows us to apply Letac's principle \cite{Letac}, which says that if $(f_i)_{i\in \nn}$ is an i.i.d. sequence of random continuous functions on some space $E$ and $\ilim n f_1\circ f_2\circ \cdots\circ f_n(x)$ almost surely converges to a constant function on $E$, then the distribution of its limit is the unique distribution stationary under iteration by an independent copy of $f_i$. In this case, taking $E=[0,e]$ and $f_i(x)=X_i^x$, we have that $\ilim n f_1\circ f_2\circ \cdots\circ f_n(x)$ converges to $T$, which is independent of $x$, and thus the distribution of $T$ is the unique one such that $X^T \eqindist T$, hence $T\eqindist Y$, so the distribution of $X$ is the inverse tower distribution of $Y$.
\end{proof} 

\begin{lemma}\label{lem:inverse_bounds}
	Suppose that $X,Y$ are independent, nonnegative, bounded random variables such that $X\eqindist Y^X$. Let $a = \inf(\supp X)$ and $b = \sup(\supp (X))$.
	
	If $b<\infty$, and $(a,b)\ne (0,1)$, then \[
	\inf(\supp (Y)) = \max(a^{\frac1a},a^{\frac1b})
	\]
	and
	\[
	\sup(\supp (Y)) = \min(b^{\frac1a},b^{\frac1b}).
	\]
	If $a=0$ and $b=1$, then $\inf(\supp(X)) = 0$ and $\sup(\supp (Y))\le1$. If $a=0$ and $b\ne 1$, we treat the quantities $a^{\frac1a}$ and $b^{\frac1a}$ as being the limit as $a$ approaches $0$ from above.
\end{lemma}
\begin{proof}
	Start by assuming $(a,b)\ne(0,1)$. Choose $\epsilon\in(0,b)$. Then note that $P(X\le a+\epsilon)$ and $P(X\ge b-\epsilon)$ are both nonnegative. Thus, we observe: \begin{eqnarray*}
		&&0 = P(X<a) = P(Y^X < a) \\&\ge& P(Y<a^{1/X}\vert X\le a+\epsilon)P(X\le a+\epsilon) +P(Y<a^{1/X}\vert X\ge b-\epsilon)P(X\ge b-\epsilon)\\
        &\ge& P\left(Y < \min(a^\frac1a,a^\frac1{a+\epsilon})\right) P(X\le a+\epsilon) + P\left(Y<\min(a^\frac1b,a^\frac1{b-\epsilon})\right)P(X\ge b-\epsilon)
	\end{eqnarray*}
        thus $P(Y<\min(a^\frac1b,a^\frac1{b-\epsilon})) =P(Y < \min(a^\frac1a,a^\frac1{a+\epsilon}))=0 $. We also have \begin{eqnarray}
		&&\nonumber P(Y < \max(({a+\epsilon})^{\frac 1{a}}, {(a+\epsilon)}^{\frac1{b}})) =P(Y^X < \max({(a+\epsilon)}^{\frac {X}{a}}, {(a+\epsilon)}^{\frac { X}b}))\\&&\ge P(Y^X < a+\epsilon) = P(X <a+\epsilon)>0.\nonumber
    \end{eqnarray}
    Therefore \[
    \max\left(\min(a^\frac1b,a^\frac1{b-\epsilon}),\min(a^\frac1a,a^\frac1{a+\epsilon})\right)\le \inf(\supp(Y)) \le \max(({a+\epsilon})^{\frac 1{a}}, {(a+\epsilon)}^{\frac1{b}})
    \]
    Letting $\epsilon$ go to $0$ we obtain the desired result for $\inf(\supp(Y))$. The proof for the supremum is similar.
    
	If $a=0$ and $b=1$, we have $\inf(\supp(Y)) = 0$ by similar reasoning to the above. We must have $\sup(\supp(Y))\le1$ because \[
	1=P(X\le 1) =P(Y^X\le 1)=P(Y\le 1)
	\]
	but nothing else can be said about $\sup(\supp(Y))$. For example, if $Y\sim\mathrm{Unif}(0,b)$ for $b\in(0,1)$ and we define $X$ to be a power tower of i.i.d. copies of $Y$ (which converges by Theorem \ref{th:main}), we have $\supp(X) = (0,1)$. 
\end{proof}
An interesting consequence of Lemma \ref{lem:inverse_bounds} is that $a^{\frac1a} \le b^{\frac1b}$ and $a^a \le b^b$ are necessary conditions for $X$ having a tower distribution. These are not guaranteed by $a\le b$. For instance, if $X$ has support $[\frac15,\frac12]$, then there can be no $Y$ such that $Y^X \eqindist X$ because $(\frac15)^{\frac15} >(\frac12)^{\frac12}$. Another consequence is that if $X$ is bounded, then $Y$ is bounded by $e^{\frac1e}$, which is the maximum of $x^{\frac1x}$ for $x\in\rr^+$.

\emph{Proof of Theorem \ref{th:inverse_problem}:} Let $T = U^r$. We observe that $P(T \le x) = \frac{x^{\frac1r}-\alpha}{\beta-\alpha}$ for $x\in(\alpha^r,\beta^r)$ and hence the probability density function of $T$ has the form \begin{equation}\label{eq:rearranged}
	p_T(x) = \begin{cases}Cx^{p-1} & x\in(a,b)\\ 0 &\text{else}\end{cases}
\end{equation}
with $p=\frac1r$, $a = \alpha^{r}$ and $b = \beta^{r}$ (or the other way around, if $r<0$), and $C$ chosen such that the total integral is $1$. Thus, it suffices to show that a random variable $T$ with a density function of the form given in (\ref{eq:rearranged}) has a tower distribution if and only if $1\in[a,b]$ and $\frac1p \in[0,\frac1{1+p\log b}]$, which is equivalent to $b <e$ and $p\ge \frac1{1-\log b}$. In that case, we claim \[
F(x) = \begin{cases}\frac{x^{bp}}{b^p}\left(1-b\log x\right)-\frac{a^p}{b^p} + \frac{a^p}{b^p} F(x^{\frac b a}) & x\in[a^{\frac1b},b^{\frac1b})\setminus\{1\}\\
	0 & x < a^{\frac1b}\\
	1 & x \ge b^{\frac1b}\\
	\frac{1-a^p}{b^p-a^p} & x = 1\end{cases}
\]
is the inverse tower distribution function of $T$. By Lemma \ref{lem:Letac}, it suffices to show that if $A$ is a random variable, independent of $T$, with distribution function $F$, then $A^T \eqindist T$. 

First, we show that if $b >e$, then $T$ does not have the tower property: By Lemma \ref{lem:inverse_bounds}, if $b\in(e,\infty)$ and $A^T\eqindist T$, then $\sup(\supp(A)) = b^{\frac1b} \le e^{\frac1e}$, which by Theorem \ref{th:main} implies that the infinite power tower formed by an i.i.d. sequence of copies of $A$ converges, but we can also see that the limit of this sequence is bounded above by $e$, so its distribution cannot be $T$. For the case of $b=\infty$, in order for the power tower of $A_i$'s to converge to $T$, we would need $\sup(\supp(A_i)) > e^{\frac1e}$. If $c>e^{\frac1e}$ and $P(A_i > c)>0$, we have $\ee (T^{\epsilon}) =\infty$ for all $\epsilon>0$ because \begin{eqnarray}
	\ee (T^{\epsilon} )&\ge& \sup_{k\in\nn}\left((c\star k)^{\epsilon} P(T > c\star k)\right) \ge \sup_{k\in\nn}\left((c\star k)^{\epsilon} P(A_i > c \mid i\in\{1,\dots,k+1\})\right) \nonumber\\&=& \sup_{k\in\nn}\left((c\star k)^{\epsilon} P(A_1>c)^{k+1}\right) = \infty.\nonumber
\end{eqnarray}
Thus, $T$ cannot have a power-law tail, because such a distribution would have $\ee( T^\epsilon) < \infty$ for sufficiently small $\epsilon$. 

If $b\le e$, then we check whether $T$ has a tower distribution. Let $G(x) = P(T \le x)$ be the distribution function of $T$. We note that, for $x\ne 1$ \begin{eqnarray}
	P(A^T \le x) &=& \int_\rr P(A^t \le x)p_T(t) dt = C\int_a^b P(A \le x^{\frac1t}) t^{p-1} dt\nonumber\\
	&=& C\int_a^b F(x^{\frac1t}) t^{p-1} dt = -C \int_{x^{\frac 1a}}^{x^{\frac1b}} F(u) \frac{(\log x)^p}{u(\log u)^{p+1}} du.\nonumber
\end{eqnarray}
In the final step, we make the substitution $t=\frac{\log x}{\log u}$. If $a=0$, we take $x^{\frac 1 a} = \lim\limits_{a\rightarrow0^+} x^{\frac 1 a}$. As we have seen, it suffices to check whether $A^T \eqindist T$. Thus $A$ has the inverse tower distribution of $T$ if and only if its distribution function $F(x)$ satisfies \begin{equation}\label{eq:inverse_equation}
	G(x) = -C (\log x)^{p}\int_{x^{\frac 1 a}}^{x^{\frac1b}} \frac{F(u) }{u(\log u)^{p+1}} du\,,
\end{equation}
so we begin by looking for solutions to this equation. It is not hard to show that there is no distribution function satisfying the equation for $p=0$, which in turn implies that $\exp(U)$ never has a tower distribution for uniform $U$. However, this is not necessary for our theorem, and is left to the interested reader to check.

Assuming $p\ne 0$, we observe that (\ref{eq:inverse_equation}) implies \begin{equation}\label{eq:inverse_equation_implication}
	F(x) = \frac{p}{Cb^p} G(x^b) - \frac{ x^b \log x}{Cb^{p-1}} G'(x^b) + \frac{a^p}{b^p}F(x^{\frac ba})
\end{equation}
almost everywhere, because
\begin{alignat}4
	\frac{G(x)}{(\log x)^p} &= -C \int_{x^{\frac1a}}^{x^{\frac1b}} \frac{F(u)}{u(\log u)^{p+1}} du\nonumber\\
	\dd x \left(\frac{G(x)}{(\log x)^p}\right) &= -C\dd x \int_{x^{\frac1a}}^{x^{\frac1b}} \frac{F(u)}{u(\log u)^{p+1}} du\nonumber\\
	\frac1{(\log x)^p}G'(x)-\frac{p}{x(\log x)^{p+1}}G(x) &= - C  \frac {x^{\frac1b - 1}}b\frac{ F(x^{\frac1b})}{\frac{x^{\frac1b}}{b^{p+1}}(\log x)^{p+1}} + C\frac{x^{\frac1a - 1}}a \frac{ F(x^{\frac1a})}{\frac{x^{\frac1a}}{a^{p+1}}(\log x)^{p+1}}\nonumber\\
	G'(x)-\frac{p}{x\log x}G(x) &=- \frac{C b^p}{x \log x} F(x^{\frac1b}) + \frac{C a^p}{x \log x} F(x^{\frac1a})\nonumber\\
	F(x^{\frac1b}) &= \frac{p}{C b^p} G(x) - \frac{x \log x}{Cb^p} G'(x) + \frac{a^p}{b^p}F(x^{\frac1a}).\nonumber
\end{alignat}
We do not necessarily have (\ref{eq:inverse_equation_implication}) implies (\ref{eq:inverse_equation}) because of the differentiation step. We do, however, have that (\ref{eq:inverse_equation_implication}) implies \[
-C \int_{x^{\frac1a}}^{x^{\frac1b}} \frac{F(u)}{u(\log u)^{p+1}} du = K + \frac{G(x)}{(\log x)^p}.
\]
We will start by showing that if $p>0$ and $F(x)$ is a distribution function satisfying (\ref{eq:inverse_equation_implication}), then $K$ in the above equation is $0$, hence (\ref{eq:inverse_equation}) is satisfied by $F$, and we conclude that $F$ is the distribution function for the inverse tower distribution of $T$. Observe that $p>0$ implies $\ilim x \frac{G(x)}{(\log x)^p} = 0$. If $F$ is distribution function, then in particular $|F(x)|\le1$, so we have\begin{eqnarray}
	|K| &=& \ilim x\abs{-C \int_{x^{\frac1a}}^{x^{\frac1b}} \frac{F(u)}{u(\log u)^{p+1}} du} \nonumber\\
	&\le& |C| \ilim x \int_{x^{\frac1b}}^{x^{\frac1a}} \frac1{u (\log u)^{p+1}} du \nonumber\\
	&=& |C| \ilim x \frac{1}{-p}\left(\frac{1}{(\log(x^{\frac1a}))^p} - \frac{1}{(\log(x^{\frac1b}))^p}\right) = |C|\ilim x \frac{b^p-a^p}{p(\log x)^p}\nonumber\\
	&=&0.\nonumber
\end{eqnarray}
Hence, to prove the theorem it suffices to show that there exists a distribution function satisfying (\ref{eq:inverse_equation_implication}) if and only if $0\le a \le b < e$ and $p\ge \frac1{1-\log b}$. 

We first show this for $a=0$, then deal with the case of $a>0$. We note that when $a=0$, we must have $p>0$, because otherwise $x^{p-1}$ does not have a finite integral. Then, (\ref{eq:inverse_equation_implication}) gives us $F$ in closed form \[
F(x) = \begin{cases}
	\frac{x^{bp}}{b^p}(1-b\log x) & x< b^{\frac1b}\\
	1 & x \ge b^{\frac1b}\end{cases}
\] 
which is the form claimed in the theorem, but is not necessarily a distribution function. If $b<1$, we observe that \[
\lim\limits_{x\rightarrow (b^{\frac1b})^-} F(x) = \frac{b^{\frac1b \cdot bp}}{b^p}(1-b\log b^{\frac1b}) = 1 - \log b > 1
\]
which means that $F(x)>1$ just below $b^{\frac1b}$, and therefore $F$ cannot be a distribution function, so for $a=0$ and $b<1$, $T$ does not have a tower distribution. If $b\ge 1$, we do have $F(x)\in[0,1]$ for all $x$, so it suffices to check whether $F$ is non-decreasing. Looking at the derivative, we can see that for $x\in(0,b^{\frac1b})$ \[
F'(x) = \frac{p x^{bp-1}}{b^{p-1}}(1-b\log x) - \frac{x^{bp-1}}{b^{p-1}} = \frac{x^{bp-1}}{b^{p-1}}(p-1-pb\log x)
\]
so $F$ is non-decreasing if and only if $p-1-pb\log x\ge 0$ for all $x\in(0,b^{\frac1b})$. Since that is a decreasing function, it suffices to check at $x=b^{\frac1b}$, hence $F$ is a distribution function if and only if \[
p-1-p\log b \ge 0 \hspace{6pt}\Longleftrightarrow \hspace{6pt} p\ge \frac1{1-\log b}
\]
which is the claim of the theorem.
Next, we deal with the case of $a>0$ and $1\in[a,b]$, after which we will conclude by showing that when $a>0$ and $1\notin[a,b]$, equation (\ref{eq:inverse_equation_implication}) has no solution, completing the proof of the theorem.

Similarly to the case of $a=0$, when $0<a\le 1\le b\le e$ we must check that $F(x)$ as defined by (\ref{eq:inverse_equation_implication}) is a distribution function. Since $b\ge 1$ and $a\le 1$, Lemma \ref{lem:inverse_bounds} implies that $F(x) = 0$ for $x<a^{\frac1b}$ and $F(x) =1$ for $x>b^{\frac1b}$, so we have \[
F(x)=\begin{cases}
	\frac{x^{bp}}{b^p}(1-b\log x) - \frac{a^p}{b^p} + \frac{a^p}{b^p}F(x^{\frac b a}) & x\in[a^{\frac1b},b^{\frac1b})\setminus\{1\}\\
	0 &x < a^{\frac1b}\\
	1 & x\ge b^{\frac1b}\\
	\frac{1-a^p}{b^p-a^p} & x = 1
\end{cases},
\]
as claimed. This obviously satisfies (\ref{eq:inverse_equation_implication}), so it remains to show that it is a distribution function if and only if $p\ge \frac1{1-\log b}$. Since $F(0)=0$ and $\ilim x F(x) = 1$, we only need to check whether $F(x)$ is non-decreasing. Observe that \[\dd x x^{bp}(1-b\log x) = bx^{bp-1}(p-1-pb\log x) = bx^{bp-1}(-1 + p(1-b\log x)),\]
which implies \[
F'(x) = \begin{cases}
	bx^{bp-1}(-1 + p(1-b\log x)) + \frac{a^{p+1}}{b^{p+1}} x^{\frac b a - 1} F'(x^{\frac b a}) & x \in (a^{\frac1b}, b^{\frac1b})\setminus\{1\}\\
	0 & x \notin[a^{\frac1b},b^{\frac1b}]\\
	\text{undefined} & \text{else}
\end{cases}.
\] If $p<0$, then $(-1 + p(1-b\log x))<0$ for all $x<b^{\frac1b}$, and hence $F$ will also have negative derivative, hence we may proceed under the assumption that $p>0$. If $p\ge \frac1{1-\log b}$, then we note that $F'(x) \ge 0$ for all $x$ where the derivative exists. Conversely, if $p< \frac1{1-\log b}$, we can show that $F'(x)<0$ for some $x$. Suppose $b>1$. For $x\in(b^{\frac a{b^2}}, b^{\frac1b})$, we have \[
F(x) = \frac{x^{bp}}{b^p}(1-b\log x),
\]
which is decreasing in a neighborhood of $b^{\frac1b}$ if $p<\frac1{1-\log b}$. The case of $b=1$ and $p< \frac1{1-\log b} = 1$ is a little more difficult. Define \[
M_n = \inf\{F'(x) : F'(x) \text{ exists and }x \le a^{a^n}\}.
\]
Obviously, if $M_n<0$ for any $n$, then $F'(x)<0$ for some $x\in(0,1)$. We claim that $p<1$ implies $\ilim n M_n<0$. Assume that $M_n\ge 0$ for all $n$. Then $\ilim n M_n$ exists (since $M_n$ is non-increasing), and we have \begin{eqnarray}
	M_{n+1} &=& \inf\left( {M_n} \cup \{F'(x) \hspace{2pt}\vert\hspace{2pt} x \in (a^{{a^n}},a^{a^{n+1}}]\}\right) \nonumber\\
	&=& \inf\left( {M_n} \cup \left\{x^{p-1}(-1 + p(1-\log x)) + a^{p+1} x^{\frac 1 a - 1} F'(x^{\frac 1 a}) \hspace{2pt}\vert\hspace{2pt} x \in (a^{a^n},a^{a^{n+1}}]\right\}\right)\nonumber\\
	&\le& \sup\left\{x^{p-1}(-1 + p(1-\log x)) \hspace{2pt}\vert\hspace{2pt}  x \in (a^{a^n},a^{a^{n+1}}]\right\} + \left(a^{p+1} (a^{a^{n+1}})^{\frac 1 a - 1} \right)M_n\nonumber
\end{eqnarray}
and therefore \begin{eqnarray}
	&&\nonumber\ilim n  M_n \\&\le& \ilim n\left(\sup\left\{x^{p-1}(-1 + p(1-\log x)) \hspace{2pt}\vert\hspace{2pt}  x \in (a^{a^n},a^{a^{n+1}}]\right\} + \left(a^{p+1} (a^{a^{n+1}})^{\frac 1 a - 1} \right)M_n\right)\nonumber\\ &=&-1 + p + a^{p+1} \ilim n M_n\nonumber
\end{eqnarray}
which can be rearranged to obtain \[
\ilim n M_n \le \frac{-1+p}{1-a^{p+1}}
\]
which is negative if $p< 1 = \frac1{1-\log b}$, a contradiction. Therefore $M_n<0$ for some $n$, which in implies that there exists $x\in (0,1)$ such that $F'(x) < 0$. Hence we have that $F$ is not a distribution if $p<1$, as claimed.

Thus, it remains to see what happens at the points where $F$ is not differentiable. These points are all those of the form $a^{\frac{a^n}{b^{n+1}}}$ or $b^{\frac{a^n}{b^{n+1}}}$ for $n\in \nn$. Define \[
F_\pm(x) = \lim\limits_{h\rightarrow0^\pm} F(x+h)
\]
to be the left- and right-limits of $F$. By induction, we can show \[
F_+(a^{\frac{a^n}{b^{n+1}}}) \ge F_-(a^{\frac{a^n}{b^{n+1}}}) 
\]
for all $n$. This is trivial for $n=0$ because \[
F_+(a^{\frac1b}) = \frac{a^p}{b^p}(1-\log a) - \frac{a^p}{b^p} = -\frac {a^p}{b^p}\log a \ge 0 = F_-(a^{\frac1b})
\]
and the recurrence relation for $F$ gives us the inductive step: \begin{eqnarray}
	F_+(a^{\frac{a^n}{b^{n+1}}}) &\ge& F_-(a^{\frac{a^n}{b^{n+1}}}) \nonumber\\
	\frac{a^{p\frac{a^n}{b^n}}}{b^p} - \frac{a^p}{b^p} + \frac{a^p}{b^p}F_+(a^{\frac{a^n}{b^{n+1}}}) &\ge& \frac{a^{p\frac{a^n}{b^n}}}{b^p} - \frac{a^p}{b^p} + \frac{a^p}{b^p}F_-(a^{\frac{a^n}{b^{n+1}}})\nonumber\\
	F_+(a^{\frac{a^{n+1}}{b^{n+2}}}) &\ge& F_-(a^{\frac{a^{n+1}}{b^{n+2}}}).\nonumber
\end{eqnarray}
We similarly have \[
F_+(b^{\frac{a^n}{b^{n+1}}}) \ge F_-(b^{\frac{a^n}{b^{n+1}}}) 
\]
for all $n$, and hence $F$ is non-decreasing in all the cases where $p\ge \frac1{1-\log b}$, as claimed. 

Finally, we check that (\ref{eq:inverse_equation_implication}) has no solutions when $1\notin[a,b]$ and $a>0$. To do this, we notice that for $x\in(a^{1/b},b^{1/b})$, (\ref{eq:inverse_equation_implication}) becomes \begin{equation}\label{eq:recursive-form-1}
    F(x) = \frac{x^{bp}}{b^p} - \frac{a^p}{b^p} - \frac{b x^{bp}\log(x)}{b^p} + \frac{a^p}{b^p}F(x^{\frac b a})
\end{equation}
or equivalently, when $x\in(a^{1/a},b^{1/a})$ \begin{equation}\label{eq:recursive-form-2}
    F(x) = -\frac{x^{ap}}{a^p} + 1 + \frac{a x^{ap}\log(x)}{a^p} + \frac{b^p}{a^p}F(x^{\frac a b}).
\end{equation}
Suppose $a>1$. Define $f_n$ recursively by $f_0(x) = 1$ and \begin{equation*}
    f_n(x)= \frac{x^{bp}}{b^p} - \frac{a^p}{b^p} - \frac{b x^{bp}\log(x)}{b^p} + \frac{a^p}{b^p}f_{n-1}(x^{\frac b a}).
\end{equation*}
Note that, if $p>0$, $\ilim x f_n(x) = -\infty$ for all $x$, and if $p<0$, $\lim\limits_{x\rightarrow 0^+} f_n(x) = \infty$. (These may be shown inductively.)

By Lemma \ref{lem:inverse_bounds}, we have that $F(x) = 1$ for $x>b^{\frac1b}$. Therefore, for $x\in (b^{\frac a{b^2}},b^{\frac1b})\cap(a^{\frac1b},b^{\frac1b})$, equation (\ref{eq:recursive-form-1}) implies \[
F(x)=\frac{x^{bp}}{b^p} - \frac{a^p}{b^p} - \frac{b x^{bp}\log(x)}{b^p} + \frac{a^p}{b^p}F(x^{\frac b a}) = \frac{x^{bp}}{b^p} - \frac{a^p}{b^p} - \frac{b x^{bp}\log(x)}{b^p} + \frac{a^p}{b^p}=f_1(x)
\]
and inductively $F(x) = f_n(x)$, for all $x\in (b^{\frac {a^n}{b^{n+1}}},b^{\frac{a^{n-1}}{b^{n}}})\cap(a^{\frac1b},b^{\frac1b})$. We observe that there exists $N\in \nn$ such that $I=(b^{\frac {a^N}{b^{N+1}}},b^{\frac{a^{N-1}}{b^{N}}}) \cap (a^{\frac1a}, a^{\frac b {a^2}}) \ne \emptyset$ because  $\bigcup\limits_{n\ge 1}(b^{\frac {a^n}{b^{n+1}}},b^{\frac{a^{n-1}}{b^{n}}})$ is open and dense in $(1,b^{\frac1b})$.

By Lemma \ref{lem:inverse_bounds}, $F(x) =0$ for $x<a^\frac1a$, so by equation (\ref{eq:recursive-form-2}), for all $x\in I$ \[
f_N(x) = F(x) = -\frac{x^{ap}}{a^p} + 1 + \frac{a x^{ap}\log(x)}{a^p} + \frac{b^p}{a^p}F(x^{\frac a b}) = -\frac{x^{ap}}{a^p} + 1 + \frac{a x^{ap}\log(x)}{a^p}.
\]
However, this is a contradiction because $f_N(x)$ and $-\frac{x^{ap}}{a^p} + 1 + \frac{a x^{ap}\log(x)}{a^p}$ are distinct analytic functions on $(0,\infty)$, so they cannot be equal on any open set. That they are distinct can be seen by looking at the limiting behavior as $x$ goes to $0$ or goes to $\infty$. When $p>0$, the limit as $x$ goes to infinity of $f_N(x)$ is $-\infty$, while the limit of $-\frac{x^{ap}}{a^p} + 1 + \frac{a x^{ap}\log(x)}{a^p}$ is $\infty$. Similarly, when $p<0$, the limits as $x$ goes to $0$ of the two functions are unequal. Thus, when $a>1$, we have that $T$ does not have the tower property. The proof for $b<1$ and $a>0$ is similar.

\qed

\section{Open questions}
There remain many open questions about random power towers. As far as convergence is concerned, Theorem \ref{th:main} is quite broad. It only leaves out the case of $a=1$ and $b=\infty$. While our condition for convergence when $a=1$ and $b>e^{\frac1e}$ does not involve $b$ in any way, our proof required boundedness of $A_1$. However, there do exist unbounded distributions for $A_1$ such that the corresponding power tower converges.
\begin{example}
	Let $\{A_i\}_{i\in\nn}$ be an i.i.d. sequence with distribution given by \[
	A_1 = \begin{cases}
		e\star n & \mathrm{ w.p. }\hspace{3pt}\frac1{2^{n+1}}\hspace{2pt}\text{for }n\in\nn\\
		e^{\frac1{e\star (16 n)}}&\mathrm{ w.p. }\hspace{3pt}\frac1{2^{n+1}}\hspace{2pt}\text{for }n\in\nn
	\end{cases}
	\]
	and let $T_i=\et\limits_{k=1}^i A_k$ be the corresponding power tower. Then $T_i$ converges almost surely.
\end{example}
To show convergence of $T_i$, it suffices to show that $\ilim i \ee \log^\star(T_i) < \infty$, since this implies that the probability of $T_i$ going to infinity is $0$. Note that \[
\ee \log^\star(A_1) = \sum_{n=1}^\infty \frac{n}{2^{n+1}} = 1 < \frac{8(1+\sqrt{1-\log 2})}{\log 2} \approx 17.94...
\]
and we show inductively that $\ee \log^\star T_i < \frac{8(1+\sqrt{1-\log 2})}{\log 2}$. Assuming this is true for $i$, we prove it for $i+1$ (making use of the fact that $\log^\star(x+y) \le 1 + \log^\star x + \log^\star y$ in the third line):\begin{eqnarray}
	&&\ee \log^\star(T_{i+1}) = \sum_{n=1}^\infty\frac1{2^{n+1}} \ee\log^\star\left((e \star n)^{T_i}\right) + \sum_{n=1}^\infty\frac1{2^{n+1}} \ee\log^\star\left(e^{T_i/(e\star (16 n))}\right) \nonumber\\
	&=& 1 + \sum_{n=1}^\infty\frac1{2^{n+1}} \ee\log^\star\left((e \star (n-2) ) + \log T_i\right) + \frac12 + \sum_{n=1}^\infty\frac1{2^{n+1}} \ee\log^\star\left(\frac{T_i}{e\star (16 n)}\right)\nonumber\\
	&\le& \frac32 + \sum_{n=1}^\infty\frac1{2^{n+1}}\left(n-2 + \ee \log^\star(T_i)\right) + \sum_{n=1}^\infty \frac1{2^{n+1}} \ee \log^\star\left(\frac{T_i}{e\star (16 n)}\right)\nonumber\\
	&\le& \frac32 + \sum_{n=1}^\infty\frac1{2^{n+1}}\left(n-2 + \ee \log^\star(T_i)\right) + \sum_{n=1}^\infty \frac1{2^{n+1}} \left(1+P(T_1\ge e\star (16 n))\ee \log^\star(T_i)\right)\nonumber\\
	&=& 2 + \frac12\ee\log^\star T_i + \sum_{n=1}^\infty\frac{1}{2^{n+1}}P(\log^\star T_1\ge (16 n))\ee \log^\star(T_i)\nonumber\\
	&\le& 2 + \frac12\ee\log^\star T_i + \sum_{n=1}^\infty\frac{1}{2^{n+1}}\frac{(\ee \log^\star(T_i))^2}{16 n} = 2 + \frac12\ee\log^\star T_i + \frac{\log 2}{32} (\ee \log^\star T_i)^2\nonumber\\
	&<& \frac{8(1+\sqrt{1-\log 2})}{\log 2}.\nonumber
\end{eqnarray}
This, and other similar examples, suggest that in order to check convergence of $T_i$ for unbounded $A_1$, we need to have some comparison between the weight of the distribution of $A_1$ at $1$ and at infinity. I conjecture that condition 2 from Theorem \ref{th:main} may be relaxed by replacing $b<\infty$ with $\ee \log^\ast A_1 < \infty$.

We also did not touch on questions of convergence rate at all. For the "contracting on average" random power towers, Diaconis and Freedman's results also give us an exponential convergence rate. Are there cases that do not have an exponential type of convergence? If so, what convergence rates are possible?

We only scratched the surface of the inverse question. If $T$ were a mixture of powers of uniform distributions, we should be able to apply similar techniques to those used in the proof of Theorem \ref{th:inverse_problem}, though the complexity of checking monotonicity increases substantially. Can one tell in general whether $T$ has a tower distribution by approximating it with uniform distributions?

%%%%%%%%%%%%%%%%%%%%%%%%%%%%%%%%%%%%%%%%%%%%%%%%%%%%%%%%%%%%%%%%%%%
%%                                                               %%
%% Use the two commands below for producing your bibliography    %%
%% with bibtex, then comment again the commands and include the  %%
%% content of the .bbl file in this file below the commands.     %%
%%                                                               %%
%%%%%%%%%%%%%%%%%%%%%%%%%%%%%%%%%%%%%%%%%%%%%%%%%%%%%%%%%%%%%%%%%%%

%\bibliographystyle{amsplain}
%\bibliography{bibliography}

\begin{thebibliography}{99}
\bibitem{Bachman}
G.~Bachman. ``Convergence of Infinite Exponentials," {\em Pacific Journal of Mathematics}, vol.~169, no.~2, pp.~219-233, 1995.

\bibitem{ARTICLE:7}
I.N.~Baker and P.J.~Rippon, ``A note on infinite exponentials,'' {\em Fibonacci Quarter}, vol.~23, pp.~106--112, 1985.

\bibitem{ARTICLE:8}
I.N.~Baker and P.J.~Rippon, ``Iterating exponential functions with cyclic exponents,'' {\em Mathematical Proceedings of the Cambridge Philosophical Society}, vol.~105, pp.~357--375, 1988.

\bibitem{Barrow}
D.F.~Barrow, ``Infinite Exponentials,'' {\em The American Mathematical Monthly}, vol.~43, no.~3, pp.~150--160, 1936.

\bibitem{article:9}
C.~Bender and J.~Vinson, ``Summation of power series by continued exponentials,'' {\em Journal of Mathematical Physics}, vol.~37, pp.~4103--4119, 1996.

\bibitem{ARTICLE:6}
R.~M. Corless, G.~H. Gonnet, D.~E.~G. Hare, D.~J. Jeffrey, and D.~E. Knuth,
``On the lambert w function,'' {\em Advances in Computational Mathematics},
vol.~5, pp.~329--359, 1996.

\bibitem{BOOK:1}
T.~H. Cormen, C.~E. Leiserson, R.~L. Rivest, and C.~Stein, {\em Growth of
	Functions: Standard notations and common functions}, pp.~58--59.
\newblock Cambridge, Massachusetts: MIT Press, 2008.

\bibitem{ARTICLE:1}
P.~Diaconis and D.~Freedman, ``Iterated random functions,'' {\em SIAM Review},
vol.~41, no.~1, pp.~45--76, 1999.

\bibitem{Dzurina}
J.~Dzurina, ``Oscillation of second order
advanced differential equations,'' {\em Electronic Journal of Qualitative Theory of Differential Equations}, vol.~20, pp.~1--9, 2018.

\bibitem{ARTICLE:2}
G.~Eisenstein, ``Entwicklung von $\alpha^{\alpha^{\alpha^{\cdots}}}$,'' {\em
	Journal fur die Reine und Angewandte Mathematik}, vol.~28, pp.~49--52, 1844.
 
\bibitem{ARTICLE:4}
L.~Euler, ``De formulis exponentialibus replicatis,'' {\em Acta Academiae
	Scientiarum Imperialis Petropolitanae}, pp.~38--60, 1778.

\bibitem{knoebel}
A.~Knoebel. ``Exponentials reiterated," {\em The American Mathematical Monthly}, vol.~88, no.~4, pp.~235-252, 1981.

\bibitem{Letac}
G.~Letac, ``A contraction principle for certain markov chains and its
applications,'' {\em Contemporary Mathematics}, vol.~50, pp.~263--273, 1986.

\bibitem{ARTICLE:5}
D.~Steinsaltz, ``Locally contractive iterated function systems,'' {\em The
	Annals of Probability}, vol.~27, no.~4, pp.~1952--1979, 1999.

\bibitem{article:3}
A.~Tarski, ``A lattice-theoretical fixpoint theorem and its applications,''
{\em Pacific Journal of Mathematics}, vol.~5, pp.~285--309, 1955.

\bibitem{Thron}
W.~Thron, ``Convergence of infinite exponentials with complex elements,'' {\em Proceedings of the American Mathematical Society}, vol.~8, pp.~1040--1043, 1957.

\end{thebibliography}

% add below the content of your .bbl file produced by bibtex.

%%%%%%%%%%%%%%%%%%%%%%%%%%%%%%%%%%%%%%%%%%%%%%%%%%%%%%%%%%%%%%%%%%%
%%                                                               %%
%% You may add acknowledgments (optional).                       %%
%%                                                               %%
%%%%%%%%%%%%%%%%%%%%%%%%%%%%%%%%%%%%%%%%%%%%%%%%%%%%%%%%%%%%%%%%%%%
\begin{acks}
I would like to thank Laurent Saloff-Coste, Persi Diaconis, Dan Dalthorp, and Lisa Madsen for their support and helpful comments. I would especially like to thank Dan Dalthorp for his work on the figures. I would also like to thank the anonymous reviewer for many helpful suggestions, especially on how to improve the introduction.
\end{acks}

%%%%%%%%%%%%%%%%%%%%%%%%%%%%%%%%%%%%%%%%%%%%%%%%%%%%%%%%%%%%%%%%%%%
%%                                                               %%
%% You have reached the end of your document.                    %%
%%                                                               %%
%%%%%%%%%%%%%%%%%%%%%%%%%%%%%%%%%%%%%%%%%%%%%%%%%%%%%%%%%%%%%%%%%%%

\end{document}